\documentclass[amstex,amstext,amsfonts,epsf,12pt] {amsart}
\usepackage{epsfig,latexsym,amsmath,amssymb,dsfont}
\usepackage{graphicx}
\usepackage{graphicx}
\usepackage{epsfig}%
\usepackage{amsfonts}%
\usepackage{amssymb}

\usepackage{color}
\usepackage[normalem]{ulem}

\newtheorem{theorem}{Theorem}[section]

\newtheorem{remark}[theorem]{Remark}

\newcommand\be{\begin{eqnarray*}}
\newcommand\ee{\end{eqnarray*}}
\newcommand\ben{\begin{eqnarray}}
\newcommand\een{\end{eqnarray}}

\makeatletter\@addtoreset{equation}{section}\makeatother
\makeatletter\@addtoreset{figure}{section}\makeatother
\makeatletter\@addtoreset{table}{section}\makeatother
\def\be{\begin{eqnarray*}}
\def\ee{\end{eqnarray*}}
\def\ben{\begin{eqnarray}}
\def\een{\end{eqnarray}}
\def\IntO{\int\limits_\Omega}

\def\Norm#1{\left\|#1\right\|}
\def\NNN {{\boldsymbol{ |}\!\!|}}
\def\wh{\widehat}
\def\wt{\widetilde}

\def\Rd{{\mathbb R}^d}
\def\Md{{\mathbb M}^{d\times d}}
\def\cA{{\mathcal A}}
\def\cB{{\mathcal B}}
\def\cE{{\mathcal E}}

\def\cV{{\mathcal V}}
\def\abig{{\bf A}}
\def\lamin{{\lambda}_\ominus}
\def\lamax{{\lambda}_\oplus}
\def\laomin{{\lambda}^\circ_\ominus}
\def\laomax{{\lambda}^\circ_\oplus}

\def\ed{\end{document}}
\def\dvg{{\rm div}}
\def\oH{\stackrel{\circ}{H}{}\!}

\def\Frame#1{
\begin{tabular}{|c|} \hline
   \\
$\
#1
$
\\
\\ \hline
\end{tabular}
}
\def\phix{\phi^{(1)}}
\def\phiy{\phi^{(2)}}
\def\Wx{W^{(1)}}
\def\Wy{W^{(2)}}
\def\fx{f^{(1)}}
\def\fy{f^{(2)}}
\def\ax{a^{(1)}}
\def\ay{a^{(2)}}

\def\bgamma{\boldsymbol\gamma}
\def\bvarsigma{\boldsymbol\varsigma}
\def\bg{\boldsymbol g}
\begin{document}

\title[Rank structured method]{Rank structured approximation method for 
quasi--periodic elliptic problems}

\author{B. Khoromskij}
\address{Max Planck Institute for Mathematics in the Sciences, 
Inselstr. 22-26, 04103, Leipzig, Germany; \textit{E-mail:} bokh@mis.mpg.de}
\author{S. Repin}
\address{V.A. Steklov Institute of Mathematics, Fontanka 27,
191 011 St. Petersburg, Russia, and University of Jyv\"askyl\"a, Finland ; 
\textit{E-mail:} repin@pdmi.ras.ru; serepin@jyu.fi}
 
%\date{}

\maketitle
\begin{abstract}
We consider an iteration method for solving an elliptic type boundary value problem 
$\mathcal{A} u=f$, where
 a positive definite operator $\mathcal{A}$ is generated by a quasi--periodic structure
with rapidly changing coefficients (typical period is characterized by a small
parameter $\epsilon$) . The method is based on using a simpler operator $\mathcal{A}_0$ 
(inversion of $\mathcal{A}_0$ is much simpler than inversion of $\mathcal{A}$), 
which can be viewed as a preconditioner for $\mathcal{A}$. We prove contraction of the 
iteration method and establish explicit estimates of the contraction
factor $q$. Certainly the value of $q$ depends on the difference between $\mathcal{A}$ and 
$\mathcal{A}_0$. 
For typical quasi--periodic structures, we establish simple relations that 
suggest an optimal $\mathcal{A}_0$ (in a selected set of "simple" structures) and compute the
corresponding contraction factor. Further, this allows us to
deduce fully computable two--sided a posteriori estimates able to control
numerical solutions on any iteration.  The method is especially efficient
if the coefficients of $\mathcal{A}$ admit low rank representations and algebraic operations
are performed in tensor structured formats.
Under moderate assumptions the storage and solution complexity
of our approach depends only weakly (merely linear-logarithmically) 
on the frequency parameter $1/\epsilon$, 
providing the FEM approximation of the order of $O(\epsilon^{1+p})$, $p>0$.
\end{abstract}



\section{Introduction}

\begin{figure}[htbp]
%\centering
%\includegraphics[width=4.0cm]{periodic2-eps-converted-to.pdf}
% \includegraphics[width=3.2cm]{Coef_X_StepSin-eps-converted-to.pdf}\quad
%%%% \includegraphics[width=4.0cm]{Fig_XSignSin_L11-eps-converted-to.pdf}
% \includegraphics[width=3.4cm]{Coef_XSinSin-eps-converted-to.pdf}\;
\includegraphics[width=7.0cm]{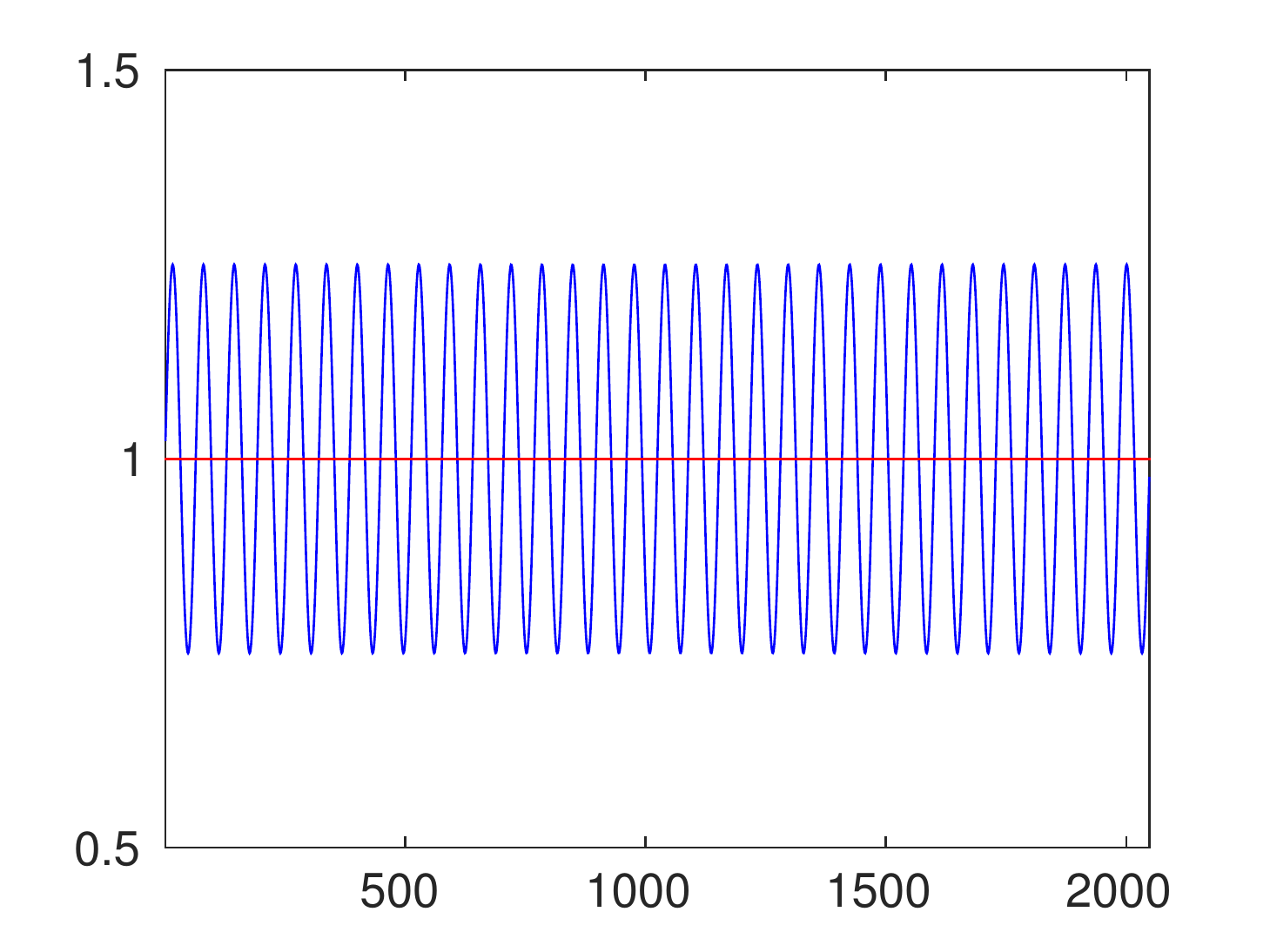}
\includegraphics[width=7.0cm]{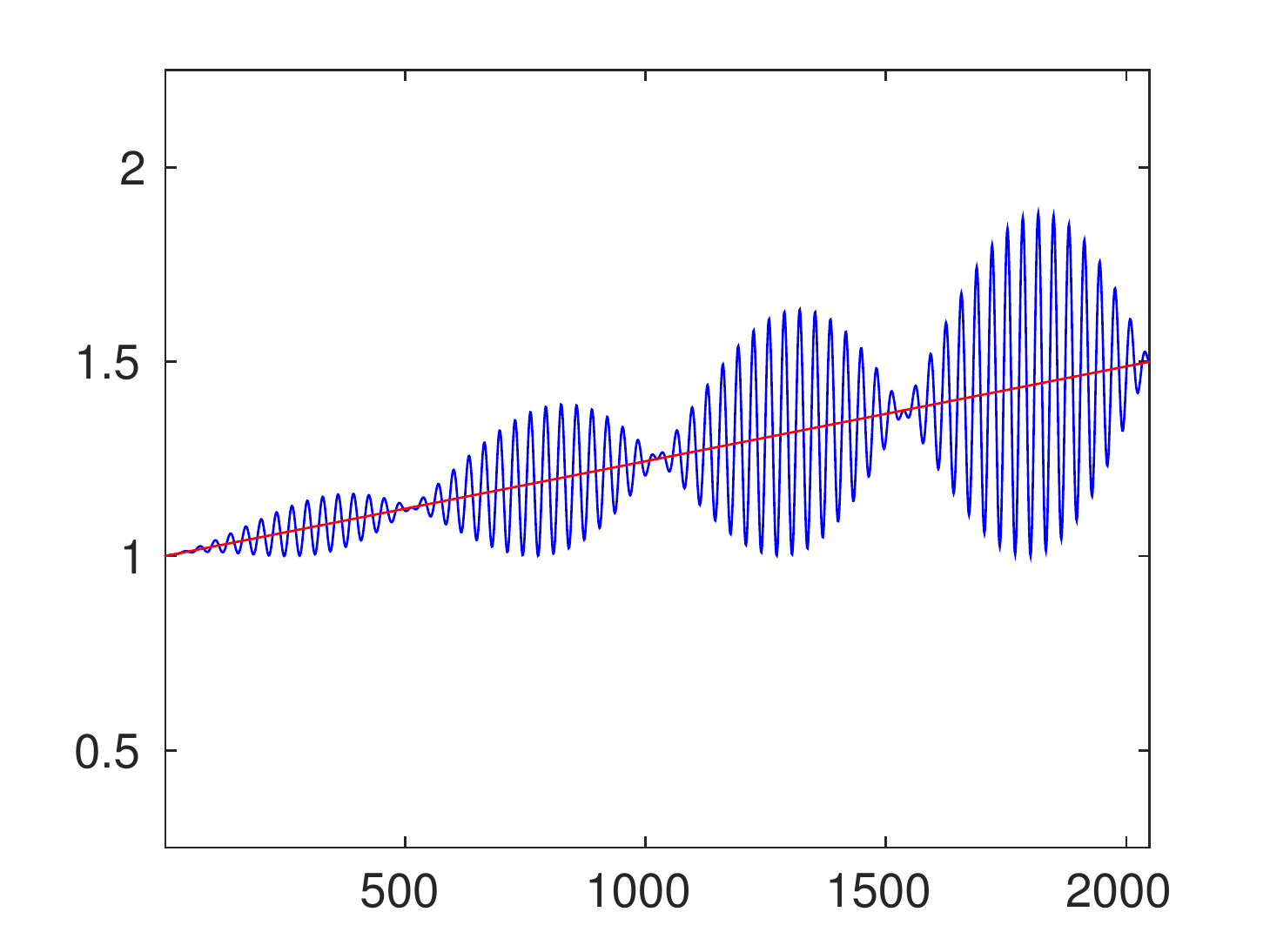}
\caption{\small Examples of periodic and modulated periodic coefficients in 1D.
}
\label{fig:1DPeriodStruct2}
\end{figure}

\begin{figure}[htbp]
%\centering
%\includegraphics[width=4.0cm]{periodic2-eps-converted-to.pdf}
% \includegraphics[width=3.2cm]{Coef_X_StepSin-eps-converted-to.pdf}\quad
%%%% \includegraphics[width=4.0cm]{Fig_XSignSin_L11-eps-converted-to.pdf}
% \includegraphics[width=3.4cm]{Coef_XSinSin-eps-converted-to.pdf}\;
%\includegraphics[width=5.4cm]{Fig_Sin_Prec_L11-eps-converted-to.pdf}
%\includegraphics[width=5.4cm]{Coef_Mod_Sin4x_X_L11-eps-converted-to.pdf}
% \includegraphics[width=4.0cm]{fig_homo-eps-converted-to.pdf}
\includegraphics[width=7.0cm]{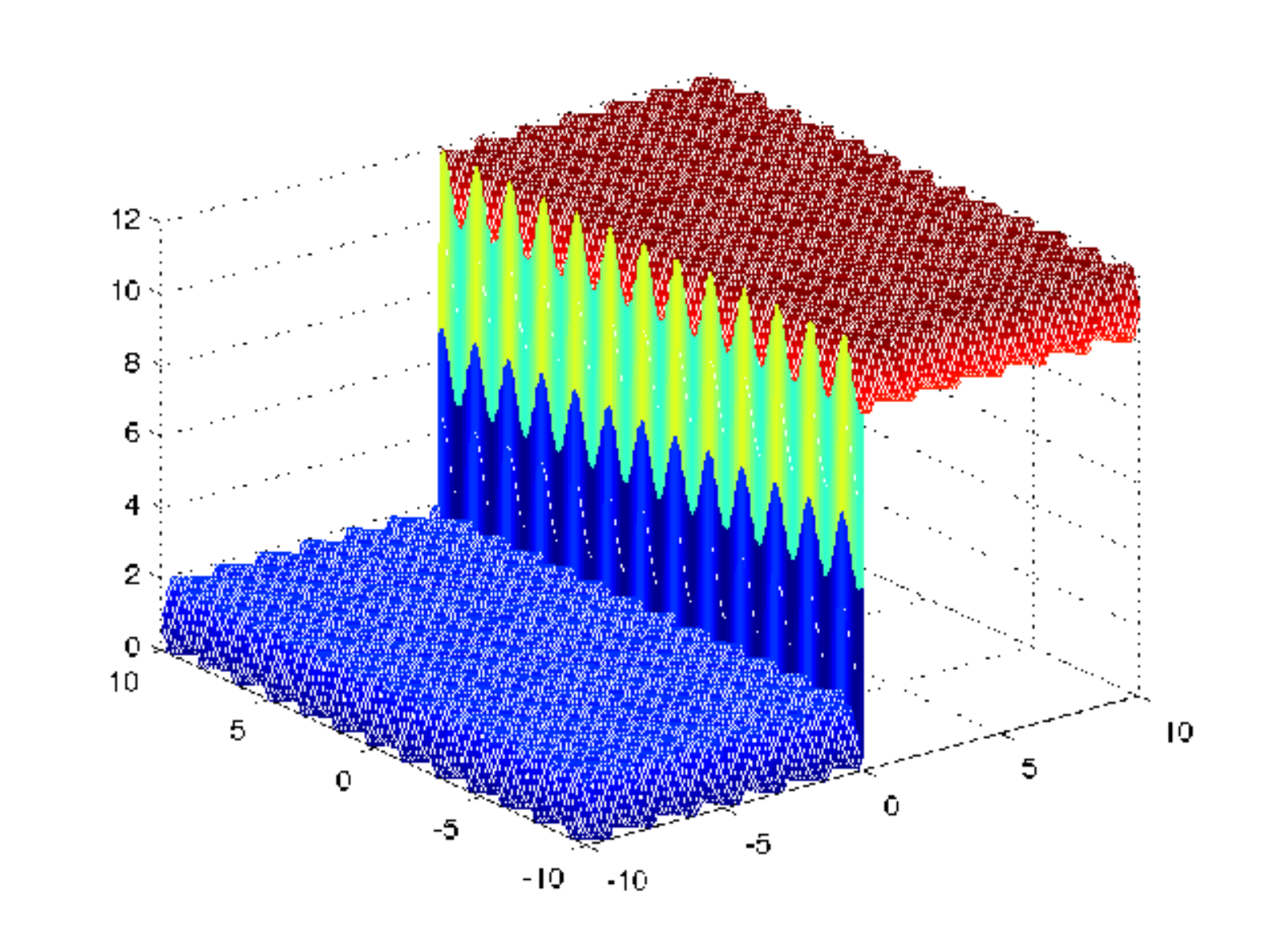}
\caption{\small An example of modulated piecewise periodic coefficients in 2D.
}
\label{fig:2DPeriodStruct2}
\end{figure}

Problems with periodic and quasi--periodic structures arise
in various natural sciences models and technical applications. 
Quantitative analysis of such problems requires 
special methods oriented towards their specific features. 
For  perfectly periodic structures, efficient
methods are developed within the framework of the homogenization theory 
(see, e.g., \cite{Bakhvalov,Bensoussan,Jikov} and other literature cited therein).
However, classical homogenization methods cover only one  class
of problems (all cells are self similar and the amount of cells is very large).
In this paper, we use a different idea and suggest another  {\em modus operandi}
for quantitative analysis 
of  boundary value problems with periodic and quasi--periodic coefficients. It
generates approximations converging (in the energy space) to the exact
solution  and provides guaranteed and computable error estimates. 
The approach is applicable to 
(see, e.g., Fig. \ref{fig:1DPeriodStruct2}, \ref{fig:2DPeriodStruct2})
\begin{enumerate}
\item
 periodic structures, in
which  the amount of cell is considerable (e.g. $10^3$--$10^4$) but not  large enough
to neglect the error generated by the respective homogenized model;
\item
 quasi--periodic structures that contain cells with defects
 and deformations;
\item
multi--periodic structures  where the coefficients reflect combined effect of several
functions with different periodicity.  
\end{enumerate}

In general terms, the idea of the method is as follows.
We consider the problem $\mathcal P$
\ben
\label{1.1}
\mathcal{A}  u=f, \quad f\in V^*,
\een
where $V$ is a reflexive Banach space with the norm $\|\cdot\|_V$, $V^*$ is the space conjugate
to $V$ (the respective duality pairing is denoted by $<v^*,v>$), 
and $\mathcal{A}:V\rightarrow V^*$ is a bounded linear operator. 
It is assumed that the operator $\mathcal{A}$ is positive 
definite and invertible, so that the problem (\ref{1.1}) is well posed. 
However, $\mathcal P$ is viewed as a very difficult problem because $\mathcal{A}$ is generated
by a complicated  physical structure, which may contain a huge amount details.  
Therefore, attempts to solve (\ref{1.1}) numerically by standard methods may lead 
to enormous expenditures. Similar difficulties arise if we wish
to  verify the quality of a numerical solution.

 Assume that the operator $\mathcal{A}$ is approximated by
a simplified positive definite operator $\mathcal{A}_\circ$ 
%(with known spectral bounds $\laomin$
%and $\laomax$) 
and the inversion of $\mathcal{A}_\circ$ is much simpler
than inversion of $\mathcal{A}$. By means of $\mathcal{A}_\circ$,
we construct an iteration method 
 based on solving a "simple" problem ${\mathcal P}_0$:
$\mathcal{A}_\circ u_\circ=g$. In other words, the method is based on the operation 
$g\rightarrow\,\mathcal{A}^{-1}_\circ g$. It also includes the operation
$v\rightarrow \,\mathcal{A}v$, which can be performed very efficiently 
by {\em tensor type decomposition methods} provided
that physical structures generated $\mathcal{A}$ have low rank
representations. We prove that iterations generate a sequence
of functions  converging to the exact solution of (\ref{1.1}) with a geometrical rate. 
Furthermore, we deduce explicitly computable and guaranteed
a posteriori error estimates adapted to this class
of problems. They evaluate the accuracy of approximations
computed on each step of the iteration algorithm.
These estimates also use only inversion of $\mathcal{A}_\circ$
and operations of the type $v\rightarrow \,\mathcal{A}v$. 
In the iteration methods and error estimates
{\em inversion of the operator $\mathcal{A}$ is avoided}.

In the paper,  we consider one  class of problems associated with
divergent type elliptic equations where
$
\mathcal{A}=Q^*\Lambda Q
$ 
and
$
\mathcal{A}=Q^*\Lambda_\circ Q.
$
Here $\Lambda: Y\,\rightarrow \, Y$ is a bounded  operator induced by a complicated quasi--periodic 
structure while $Q:V\,\rightarrow \, Y$ and $Q^*:Y\,\rightarrow \, V^*$ are conjugate
 operators, i.e.,
\ben
\label{1.2}
(y,Qw)=<Q^*y,w>\qquad \forall y\in Y\;{\rm and}\;w\in V,
\een
where $Y$ is a Hilbert space with the scalar product $(\cdot,\cdot)$ and  the norm $\|\cdot\|$. 
The operators $Q$ and $Q^*$ are induced by differential operators
or certain finite dimensional approximations of them. Henceforth, it is assumed that
$f\in \cV$, where
$\cV$ is a Hilbert space with the scalar product $(\cdot,\cdot)_\cV$. 
This space is intermediate between
$V$ and $V^*$, i.e., $V\in \cV\in V^*$.

The operator $\mathcal{A}_\circ=Q^*\Lambda_\circ Q$ contains the operator 
$\Lambda_\circ$ generated by a simplified structure.
We assume that the operators $\Lambda$ and $\Lambda_\circ$ are Hermitiam
(i.e.,
$
(\Lambda y,z)=(y,\Lambda z) 
$
and
$
(\Lambda_\circ y,z)=(y,\Lambda_\circ z)
$)
and satisfy the conditions
\ben 
\label{1.Lambda0}
&&\laomin\|y\|^2\,\leq\, (\Lambda_\circ y, y)\,\leq\,\laomax\|y\|^2\qquad
\forall y\in Y,\\
\label{1.Lambda}
&&\lamin\,\|y\|^2\;\;\leq\;\; (\Lambda y, y)\,\leq\,
\lamax\,\|y\|^2, \qquad \lamin<\lamax.
\een
 Then, the structural operators $\Lambda$ and $\Lambda_\circ$ are spectrally equivalent
\ben
\label{1.equivalence}
&&c_1 (\Lambda_\circ y, y)\;\leq\; (\Lambda y, y)\,\leq\,c_2 (\Lambda_\circ y, y),
\een
where  the constants are the minimal and maximal eigenvalues of the generalized
spectral problem $\Lambda y-\mu \Lambda_\circ y=0$.
Obviously, they satisfy the estimates
$c_1\geq\frac{\lamin}{\laomax}$ and 
$c_2\leq\frac{\lamax}{\laomin}$ (which may be rather coarse).

Concerning the operator $Q$, we assume that there exists a positive constant $c$ such that
\ben
\label{1.6}
\|Q w\|\,\geq\, c\|w\|_V\qquad\forall w\in V.
\een
%%%%%%%%%%%%%%%%%%%%%%%
Generalized solutions of the problems ${\mathcal P}$ and ${\mathcal P}_0$ are defined
by the variational identities
\ben
\label{1.7}
&&(\Lambda Q u,Qw)=<f,w>\qquad\forall \,w\in V,
\een
and
\ben
\label{1.8}
&&(\Lambda_\circ Q u_0,Qw)=<\wt f,w>\qquad\forall \,w\in V.
\een
In Sect. 2, we show that a sequence $\{u_k\}$ converging to $u$ in $V$
can be constructed by solving  problems (\ref{1.8}) with specially
constructed right hand sides $\wt f_k$ generated by the residual of (\ref{1.7}). 
In proving convergence, the key
issue is analysis of the spectral radius of the operator
  \ben
 \label{1.9}
{\mathbb B}_\rho\,:={\mathbb I}-\rho \Lambda^{-1}_\circ \Lambda,
\een
and selection of such relaxation parameter $\rho$ that provides the best convergence rate.
Moreover, iteration procedures of such a type become contracting
if the iteration parameter is properly selected. This fact
is often used in proving analytical results (e.g., see \cite{LiSt1967}, 
where classical results on existence
and uniqueness of a variational inequality has been established
by  contraction arguments ). Also, these ideas were used
in construction of various numerical methods
(see, e.g., \cite{GlLiTr}). However, achieving our goals requires more than the fact of 
contraction. We need explicit
and realistic estimates of the contraction factor (which are used in
error analysis) and a practical method of finding $\Lambda_\circ$ with minimal
$q$. The latter task leads to a special optimization problem that
defines the most efficient "simplified" operator among  a certain class
of "admissible" $\Lambda_\circ$. These questions are studied
in Sect. 3.
In general, $\Lambda$ and $\Lambda_\circ$ can be  induced by scalar,
vector, and tensors functions. We show that selection of the optimal 
structural operator $\Lambda_\circ$  is reduced to a special interpolation type problem, 
which is purely algebraical and 
does not require solving a differential problem (therefore selection of a suitable
$\Lambda_\circ$ can be done a priori). We discuss several examples
and suggest the corresponding optimal (or quasi optimal) $\Lambda_\circ$,
which guarantees convergence of the iteration sequence with {\em explicitly
known contraction factor}. 

Now, it is worth saying  about the main differences between our approach and the classical homogenization method developed for regular periodic structures. This method operates with a homogenized boundary value problem 
$Q^*\Lambda_{\rm H} Q \,u_{\rm H}=f$, where $\Lambda_{\rm H}$ is defined by means of
an auxiliary problem
with periodical boundary conditions in the cell of periodicity.   The respective solution $u_{\rm H}$
contains an irremovable (modeling) error depending on the cell diameter $\epsilon$.
Moreover, if $\epsilon$ tends to zero, then typically $u_{\rm H}$ converges to $u$ only 
weakly (e.g., in $L^2$). Getting a better convergence (e.g., in $H^1$) requires certain
corrections, which lead to other (more complicated) boundary value problems in the cell of periodicity.
The respective "corrected" solution $u^c_{\rm H}$ also contains an error. 
Typically, the error
is  proportional to $\sqrt{\epsilon}$ and can be neglected only if the amount of cells is very large.
If our method is applied to perfectly periodical structures then setting
$\Lambda_\circ:=\Lambda_{\rm H}$ is one possible option. In this case, the homogenized operator
(defined  without correction procedures) is used for a different purpose: 
{\em construction of a suitable preconditioning operator}. The latter operator generates
numerical solutions converging to the exact solution in the energy norm 
(i.e., the method is free from irremovable errors)  and can be applied for a rather wide
range of  $\epsilon$. In addition, the theory suggests other simpler ways of 
selecting suitable $\Lambda_\circ$. In this context, it is interesting to know weather or not  
the choice $\Lambda_\circ:=\Lambda_H$ always yields  minimal value of the contraction factor.
In Sect. 3, we briefly discuss this question and  present an example of that the best
$\Lambda_\circ$ 
may differ from $\Lambda_{\rm H}$ .

In Sect. 4, we deduce a posteriori estimates that
provide fully computable and guaranteed estimates of the distance to the
exact solution $u$ for any numerical approximation $u_{k,h}$ computed
for an approximation subspace $V_h$. These estimates are established by
combining functional type a posteriori estimates (see \cite{Re2000,NeRe,ReGruyter} and 
references cited therein)  and estimates generated
by the contraction property of the iteration method (see \cite{Ostrowski,Zeidler}).

The second part of the paper is devoted to a fast solution method for  the basic 
iteration problem (\ref{2.1}). The key idea consists of using 
  tensor type representations for approximations, what is quite natural if 
  both coefficients of the respective quasi--periodic structure and the right-hand side
  admit low rank tensor type representations. 
  We notice  that the amount of structures representable in terms of low rank formats is much
  larger than the amount of  periodic structures covered by the homogenization
  method. The idea of tensor type approximations of partial differential
  equations traces back to \cite{KantorovichKrylov}. In computational
  mechanics this method is known as the Kantorovich--Krylov (or extended Kantorovich) method.
  However, it is rarely used in modern numerical technologies, which are mainly
  based upon various finite element technologies. In part, this is due restrictions
  on the shape of the domain imposed by the Kantorovich method.
  Henceforth, we assume that the domain $\Omega$ satisfies these restrictions,
  i.e.,  it is a tensor type domain (e.g., rectangular) or a union of tensor type domains. 
  Certainly, this fact induces some limitations, which however could be bypassed 
  by known methods (coordinate transformation, domain decomposition, iso-geometric analysis, etc.).
  
  The recent tensor numerical methods for steady state and dynamical problems
  based on the advanced nonlinear tensor approximation algorithms have been developed in the last
  ten years.
  Literature survey on the modern tensor numerical methods for multi-dimensional PDEs
  can be found in \cite{KhorCA:09,KhorSurv:10,VeKhorTromsoe:15}.
  In the context of problems considered in the paper, 
  we are mainly concerned with another specific feature: very complicated material structure. 
  In this case, direct application of standard finite element methods suffers
  from the necessity to account huge information encompassed in coefficients (especially in
  multi dimensional problems). We show that tensor type methods allow us to
  reduce computations to a collection of one dimensional problems, 
  which can be solved very efficiently using low rank representations with the small
  storage requests. 
  Similar ideas are applied  for computing a posteriori   error estimates.
  
Section 5 discusses numerical aspects of the method and exposes several
examples. Typical behavior of  quasi-periodic coefficients is described by 
oscillation around constant, modulated oscillation around given smooth function, 
or oscillation around piecewise constant function. 

Figure \ref{fig:1DPeriodStruct2} (1D case)
represents examples of highly oscillating (left) and modulated periodic coefficients (right) functions.

Figure \ref{fig:2DPeriodStruct2} (2D case) illustrates the well separable equation 
coefficient obtained by a sum of step-type and uniformly oscillating functions.

%\color{green}
We show that specially constructed FEM type approximations of PDEs with 
slightly perturbed or regularly modulated periodic coefficients on $d$-fold 
 $\underbrace{n\times \cdots \times n}$  tensor 
grids in $\mathbb{R}^d$ may lead to the discretized algebraic
equations with the low Kronecker rank stiffness matrix of size $n^d \times n^d$, where 
$n=O(\frac{1}{\epsilon})$ is proportional to the large frequency parameter $1/\epsilon$. 
In this case the rank decomposition with respect to the $d$ spacial variables is applied,
such that the discrete solution can be calculated in the low-rank separable form,
which requires the only $O(dn)$ storage size instead of $O(n^d)=O(\frac{1}{\epsilon^d})$ 
complexity representations which are mandatory for the traditional FEM techniques 
(the latter quickly leads to the bottleneck in case of small parameter $\epsilon>0$). 

The arising linear system of equations can be solved by preconditioned iteration 
with the simple preconditioner $\Lambda_\circ$,
such that the storage and numerical costs scale almost linearly
in the univariate discrete problem size $n$, i.e., they are estimated by
\[
 O(d n  \log^p(\frac{1}{\epsilon}))\ll O(\frac{1}{\epsilon^d}),\quad p > 0,
\]
% $$
% O(\log^p(\frac{1}{\epsilon^d}))=O(d^p\log^q(\frac{1}{\epsilon})),\quad
% p\geq 1,
% $$
where $d$ is the spatial dimension.
Numerical examples in Section 5 demonstrate the stable geometric convergence
of the preconditioned CG (PCG) iteration with the preconditioner $\Lambda_\circ$ and 
confirm the the low-rank approximate separable representation to the solution
with respect to $d$ spacial variables even 
in the case of complicated quasi-periodic coefficients. 

This approach is well suited for applying the quantized-TT (QTT)
%\cred What is QTT the reader may ask? \cn
tensor approximation \cite{KhQuant:09}
to  functions discretized on large tensor grids of size proportional 
to the frequency parameter, i.e.  $n=O(1/\epsilon)$, as it was demonstrated in 
the previous paper \cite{BokhSRep:15} for the case $d=1$. 
The use of tensor-structured 
preconditioned iteration with the adaptive QTT rank truncation may lead to the 
logarithmic complexity in the grid size, $O(\log^p n)$,
see \cite{KhorCA:09,KhorSurv:10,OselDolg:12} 
for the rank-truncated iterative methods,   
\cite{VeBoKh:Ewald:14,VeKhorCorePeriod:14,DoKazKh_1DSPDE:12,KaRaSch:11} 
for various examples of the QTT tensor approximation 
to lattice structured systems, and \cite{BeKh3_Prot1:16}
for tensor approximation of complicated functions with multiple cusps in $\mathbb{R}^d$.

In Section 6, we conclude with the discussion on further perspectives of the presented approach
for 2D and 3D elliptic PDEs with quasi periodic coefficients.
%\color{black}

% \color{red}
% \noindent
% Question 1: what is PCG? Should be explained.\\
% Question 2: what means "existence of low rank solutions"? Finite dimensional problem generated by tensor approximations always has a solution if one dimensional basic functions are properly selected.\\
% Question 3: I do not like the very last sentence. First, what is "geometric
% homogenization", second we outline from the very beginning that our
% method is much more general and applicable for a wide set of quasi periodic
% structures. May be write"...for 2D 3D problems with quasi periodic coefficients?
% \color{black}

%Then,
%\ben
%\label{1.3}
%(A_0-\rho A)\zeta\cdot\zeta\,\leq\,\kappa A_\circ\zeta\cdot\zeta.
%\een
%In view of (\ref{1.3}),
%\ben
%\label{1.4}\zeta_1:=\frac{\min\{a^{(1)}_\ominus,a^{(2)}_\ominus\}
%{\max\{a^{(1)}_\ominus,a^{(2)}_\ominus\}}
%(1-\kappa)A_\circ\zeta\cdot\zeta\,\leq\,\rho A\zeta\cdot\zeta\,\leq\,\rho\|A\|\,\|\zeta\|^2.
%\een

%%%%%%%%%%%%%%%%%%%%%%%

\section{The iteration method }
Let $v\in V$ and $\rho\in {\mathbb R}_+$.
Consider the problem: find $u_v$ such that
\ben
\label {2.1}
(\Lambda_\circ Q u_v, Q w)
=\ell^\circ_v(w) - \rho\, \ell_v(w)\qquad \forall w\in V,
\een
where
\be
\ell_v(w):=(\Lambda Qv, Qw)-
<f,w>\quad{\rm and}\quad \ell^\circ_v(w):= (\Lambda_\circ Q v, Q w).
\ee
Obviously, the right hand side of (\ref{2.1}) is a bounded linear
functional on $V$, so that this problem has a unique solution $u_v$. Thus,
we have a mapping $T_\rho:V\rightarrow V$, which becomes a contraction if the parameter  $\rho$ is properly selected.
Indeed, for  any $v_1$ and $v_2$ in $V$,
 we obtain
 \ben
 \label{2.2}
 (\Lambda_\circ Q\eta, Qw)
%=(\Lambda_\circ Q\zeta, Qw) - \rho (\Lambda Q\zeta, Qw)
=(\Lambda_\circ Q\zeta-\rho \Lambda Q\zeta,Qw) )\qquad \forall w\in V,
\een
where $u_1 = T_\rho v_1$,  $u_2 = T_\rho v_2$, $\zeta:=v_1-v_2$, and $\eta:=u_1-u_2$.
Hence
\begin{multline}
\label {2.3}
\|\eta\|^2_\circ:=(\Lambda_\circ Q\eta, Q\eta)
=(\Lambda_\circ Q\zeta, Q \eta)-\rho (\Lambda Q\zeta, Q \eta)\\
=(Q\zeta, \Lambda_\circ Q\eta) -\rho (\Lambda^{-1}_\circ\Lambda Q\zeta, \Lambda_\circ Q\eta)
=(Q\zeta-\rho \Lambda^{-1}_\circ\Lambda Q\zeta, \Lambda_\circ Q\eta)\\
\leq
\|\eta\|_\circ
\left(\Lambda_\circ Q\zeta-\rho\,
\Lambda Q  \zeta, Q\zeta-\rho\,
\Lambda^{-1}_\circ\Lambda Q \zeta\right)^{1/2}.
\end{multline}
>From (\ref{2.3}) we find that
\begin{multline}
\label{2.4}
\|\eta\|^2_\circ\leq\,(\Lambda_\circ Q\zeta, Q\zeta)-2\rho (\Lambda Q\zeta, Q\zeta)+\rho^2
(\Lambda^{-1}_\circ \Lambda Q\zeta, \Lambda Q\zeta)\\
=(Q\zeta, \Lambda_\circ Q\zeta)-2\rho (\Lambda^{-1}_\circ\Lambda Q\zeta,\Lambda_\circ Q\zeta)+\rho^2
(\Lambda^{-1}_\circ\Lambda\Lambda^{-1}_\circ\Lambda Q\zeta, \Lambda_\circ Q\zeta)\\
=(({\mathbb I}-2\rho\Lambda^{-1}_\circ\Lambda+\rho^2\Lambda^{-1}_\circ\Lambda\Lambda^{-1}_\circ\Lambda)Q\zeta,\Lambda_\circ Q\zeta)=(\Lambda_\circ{\mathbb B}^2_\rho Q\zeta,Q\zeta)\\
\leq (\Lambda_\circ{\mathbb B}^2_\rho Q\zeta,{\mathbb B}^2_\rho Q\zeta)^{1/2}(\Lambda_\circ Q\zeta,Q\zeta)^{1/2},
\end{multline}
where ${\mathbb B}_\rho$ is defined by (\ref{1.2}). If   $\rho$ is selected such that
\ben
\label{2.5}
(\Lambda_\circ{\mathbb B}^2_\rho Q\zeta,Q\zeta)\,\leq\, q^2\|\zeta\|^2_\circ,
\qquad {\rm for\; some}\;q<1,
\een
then (\ref{2.4}) shows that $T_\rho$ is a contractive mapping.

It is not difficult to show that  $\rho$ satisfying  (\ref{2.5}) can be always found.
Indeed, in view of  (\ref{1.equivalence})
\begin{multline}
\label{2.6}
(\Lambda_\circ{\mathbb B}^2_\rho Q\zeta,Q\zeta)=
(\Lambda_\circ Q\zeta,Q\zeta)-2\rho(\Lambda Q\zeta,Q\zeta)+\rho^2(\Lambda\Lambda^{-1}_\circ\Lambda Q\zeta,Q\zeta)\\
\leq  (1-2\rho c_1)(\Lambda_\circ Q\zeta,Q\zeta)+\rho^2
(\Lambda^{-1}_\circ\Lambda Q\zeta, \Lambda Q\zeta).
\end{multline}
Since $\Lambda$ and $\Lambda_\circ$ are invertible with
 trivial kernels,
 $\mu$ and $y_\mu$ are an eigenvalue and the respective
eigenfunction of $\Lambda y_\mu=\mu \Lambda_\circ y_\mu$
if and only if
they are an eigenvalue and the eigenfunction
of the problem 
$\Lambda\Lambda^{-1}_\circ\Lambda y_\mu=\mu \Lambda y_\mu$.
This means that
\be
c_1(\Lambda y,y)\,\leq\,(\Lambda\Lambda^{-1}_\circ\Lambda y,y)\, \leq\,
c_2(\Lambda y,y)\,\leq\, c^2_2(\Lambda_\circ y,y).
\ee
%Since
%$
%\frac{1}{\lamax}(y,y)\leq
%(\Lambda^{-1}y,y)\,\Rightarrow\,
%(y=\Lambda\zeta)\,\Rightarrow\,
%(\Lambda\zeta, \Lambda\zeta)\,\leq\,
%\lamax 
%(\Lambda\zeta,\zeta),
%$
Hence
\be
(\Lambda^{-1}_\circ\Lambda Q\zeta, \Lambda Q\zeta)\leq 
c^2_2 \|\zeta\|^2_\circ
%\frac{1}{\laomin}\|\Lambda Q\zeta\|^2\leq
% \frac{\lamax}{\laomin}(\Lambda Q\zeta, Q\zeta)\leq
%\frac{\lamax^2}{{\laomin\!}^2}\, \|\zeta\|^2_\circ,
\ee
and (\ref{2.6}) implies
\ben
\label{2.7}
(\Lambda_\circ{\mathbb B}^2_\rho Q\zeta,Q\zeta)\,\leq\, 
\left(1-2\rho c_1 +\rho^2 c^2_2
\right)\|\zeta\|^2_\circ.
\een
Minimum of the expression in round brackets
is attained if $\rho=\rho_*:=  \frac{c_1}{c^2_2}$.
For $\rho=\rho_*$, we find that
\ben
\label{2.8}
q^2_*:=1-\frac{c^2_1}{c^2_2}\,
\leq\, \hat q^2:=1-\frac{\lamin^2{\laomin\!}^2}
{\lamax^2{\laomax\!}^2}
\in \,[0,1).
\een

 Hence, $T_\rho$ is a contractive mapping with explicitly known contraction factor $q_*$.
Well known results in the theory of fixed points (e.g., see \cite{Zeidler})
yield the following result.

\begin{theorem}
\label{th1}
For any $u_0\in V$
and $\rho=\rho^*$   the  sequence $\{u_k\}\in V$ of functions satisfying the relation
\ben
\label{2.9}
(\Lambda_\circ Q u_k, Qw)=(\Lambda_\circ Qu_{k-1}, Qw) - \rho \Bigl( (\Lambda Q u_{k-1}, Qw)-
<f,w>\Bigr) \qquad \forall w\in V
\een
converges
to $u$ in $V$ and
$\|u_k-u\|_\circ
\leq q^k_*\|u_0-u\|_\circ$
as $k\rightarrow+\infty$.
\end{theorem}
\begin{remark}
>From (\ref{2.4}) we obtain
\be
\|\eta\|^2_\circ
\leq \frac{1}{\lambda_{0,{\rm min}}}\|{\mathbb B}^2_\rho Q\zeta\|\|\zeta\|_\circ\leq 
\frac{\|{\mathbb B}^2_\rho\|}{\lambda^2_{0,{\rm min}}}\|\zeta\|^2_\circ.
\ee
This relation yields a simple (but not very sharp) estimate
of the contraction factor.
\end{remark}
For further analysis, it is convenient 
to estimate the right hand side of (\ref{2.4}) by a different method.
Let $\NNN {\mathbb B}_\rho\NNN_\circ $ denote the operator
norm
\ben
\label{2.10}
\NNN {\mathbb B}_\rho\NNN_\circ:=\sup\limits_{y\in Y}
\frac{\| {\mathbb B}_\rho y\|_\circ}{\| y\|_\circ}.
\een
Then
$
\| {\mathbb B}_\rho y\|_\circ\,\leq\,\NNN {\mathbb B}_\rho\NNN_\circ
 \|y\|_\circ
$
and
\be
(\Lambda_\circ{\mathbb B}^2_\rho y,y)\,
\leq\,\NNN {\mathbb B}_\rho\NNN^2_\circ\,\| y\|^2_\circ.
\ee
Hence, (\ref{2.4}) yields the estimate
\ben
\label{2.11}
\|\eta\|_\circ\,\leq\,
\NNN {\mathbb B}_\rho\NNN_\circ \|\zeta\|_\circ,
\een
which shows that $T_\rho$ is a contraction provided that
\ben
\label{2.12}
\NNN {\mathbb B}_\rho\NNN_\circ\,<\,1.
\een
In applications ${\mathbb B}_\rho$ is a self adjoint bounded
operator acting in a finite dimensional space, so that verification of this condition
amounts finding $\rho$ which yields the respective
spectral radius of ${\mathbb B}_\rho$ (see Section 4).

%%%%%%%%%%%%%%%%%%%%%%%%%%%%%%

\section{Selection of $\Lambda_\circ$}
In this section, we discuss how to select $\Lambda_\circ$
in order to minimize $q$ what is crucial for two major aspects
of quantitative analysis: convergence of the iteration method
and guaranteed a posteriori estimates. We assume that $V$, $\cV$, and $Y$ are spaces
of functions defined in a Lipschitz bounded domain $\Omega$
(namely $y(x)\in \mathbb T$ for a.e. $x\in \Omega$ where $\mathbb T$ may coincide
with $\mathbb R$, ${\mathbb R}^d$, or $\Md$)
and the operators $\Lambda$ and $\Lambda_\circ$ are generated by bounded
scalar functions, matrices or tensors. In this case,
\be 
(\Lambda y,y):=\IntO \Lambda(x) y\odot y\,dx,\quad
{\rm and}\quad
(\Lambda_\circ y,y):=\IntO \Lambda_\circ(x) y\odot y\,dx,
\ee
where $\odot$ denotes the respective product of scalar, vector, or tensor
functions.
In view of (\ref{2.10}) and (\ref{2.12}), the value of $\rho$ should minimize
 the quantity
$
\sup\limits_{y\in Y}
\frac{(\Lambda_\circ {\mathbb B}_\rho y,{\mathbb B}_\rho y)}{(\Lambda_\circ y,y)}.
$
This procedure yields the contraction factor
\ben
\label{3.1}
q^2={\mathcal Q}(\Lambda,\Lambda_\circ):=\inf\limits_{\rho}\,\sup\limits_{y\in Y}
\frac{\IntO \Lambda_\circ(x) {\mathbb B}_\rho(x) y\odot 
{\mathbb B}_\rho (x) y\,dx}{\IntO \Lambda_\circ(x) y\odot y\,dx},
\een
which computation is reduced to  ${\mathbb B}_\rho(x)$  to solving
algebraic problems at a.e. $x\in \Omega$, i.e.,
 \ben
 \label{3.2}
{\mathcal Q}(\Lambda,\Lambda_\circ):=\inf\limits_{\rho}\,\sup\limits_{x\in \Omega}\,\sup\limits_{\tau\in {\mathbb T}}
\frac{\Lambda_\circ(x) {\mathbb B}_\rho(x) \tau\odot 
{\mathbb B}_\rho (x) \tau}{ \Lambda_\circ(x) \tau\odot \tau}
\een

%This supremum coincides with the maximal eigenvalue of the generalized
%eigenvalue problem 
%${\mathbb B}_\rho\Lambda_\circ {\mathbb B}_\rho\,y-\mu \Lambda_\circ\,y=0$.
%% which is equivalent to
%%\be 
%(\Lambda_\circ-\rho\Lambda^{-1}_\circ\Lambda\Lambda_\circ)({\mathbb I}-\rho \Lambda^{-1}_\circ \Lambda)y-\mu\Lambda_\circ y=0.
%\ee

Let $\mathds S$ be a certain set of "simple" operators
defined a priori (e.g., it can be a finite dimensional
set formed by piece vise constant or polynomial functions).
Then, finding the best "simplified" operator amounts
solving the problem: find $\widehat\Lambda_\circ\in \mathbb S$ such that
${\mathcal Q}(\Lambda,\widehat\Lambda_\circ)$ is minimal. In other words,
optimal $\widehat\Lambda_\circ$ is defined by the problem
\begin{equation}
\label{3.3}
\Frame{\inf\limits_{\,\Lambda_0\in \mathds S,\atop \rho \in \mathbb R}\, 
\sup\limits_{ x \in \Omega  \atop \tau\in {\mathbb T}}
\frac{\displaystyle\Lambda_\circ(x) {\mathbb B}_\rho(x) \tau\odot 
{\mathbb B}_\rho (x) \tau}{\displaystyle \Lambda_\circ(x) \tau\odot \tau}\,=\,q^2.
}
\end{equation}
Notice that  (\ref{3.3}) is an algebraic
problem, which should be solved (analytically or numerically) 
before computations. The respective solution $\widehat\Lambda_\circ$
defines the best operator to be used in the iteration method (\ref{2.9}) 
and yields the respective contraction factor.
Below we discuss some  particular cases,
where analysis of this problem generates optimal (or almost
optimal) $\Lambda_\circ$.

%\subsection{Structures generated by scalar coefficients}
Problem  (\ref{3.3}) is explicitly solvable if $\Lambda_\circ$ and $\Lambda$ have a special structure, namely,
\be 
\Lambda_\circ=a_\circ(x) {\mathbb I},\qquad
\Lambda=a(x) {\mathbb I},
\ee
where ${\mathbb I}$ is the unit operator and $a_\circ(x)$ and $a(x)$ are positive bounded functions defined in  $\Omega$. Then, \be 
 {\mathbb B}_\rho(x)=(1-{\mathds h}(x)){\mathbb I},\qquad{\mathds h}(x):=\frac{a(x)}{a_\circ(x)}
\ee
 and
  \be 
\sup\limits_{\tau \in {\mathbb T}}
\frac{\left(1-\rho {\mathds h}(x)\right)^2\tau\odot\tau}
{|\tau|^2}\,=\,|1-\rho{\mathds h}(x)|^2\qquad \forall\,x\in \Omega .
\ee
Define
${\mathds h}_\ominus:=\min\limits_{x\in \Omega}{\mathds h}(x)$
and
${\mathds h}_\oplus:=\max\limits_{x\in \Omega}{\mathds h}(x)$.
It is not difficult to show that
$$
\sup\limits_{x\in \Omega} |1-\rho {\mathds h}(x)|=\max\{|1-\rho {\mathds h}_\ominus|,|1-\rho {\mathds h}_\oplus|\}.
$$
Minimization with respect to $\rho$ yields the best value
 $\rho_*=\frac{2}{{\mathds h}_\ominus+{\mathds h}_\oplus}$ and the respective value
 \ben
 \label{3.4} 
 {\mathcal Q}(\Lambda,\Lambda_\circ)=\left(\frac{{\mathds h}_\oplus-{\mathds h}_\ominus}{{\mathds h}_\oplus+{\mathds h}_\ominus}\right)^2=
 \left(\frac{1-{\mathcal J(a,a_\circ)}}
 { 1+{\mathcal J(a,a_\circ)}}\right)^2<1,\qquad {\mathcal J(a,a_\circ)}=
\frac{ {\mathds h}_\ominus}{{\mathds h}_\oplus}.
 \een
In accordance with (\ref{3.3}) identification of the optimal simplified
problem is reduced to the problem
\ben
\label{3.5}
\sup\limits_{a_0\in {\mathbb S}}\;{\mathcal J(a,a_\circ)}.
\een
where ${\mathbb S}$ is a given set of functions.

We illustrate the above relations by means of several examples. \\
{\it Example 1. Constant coefficients.}
In the simplest case, we set ${\mathbb S}=P^0$, i.e., $a_0$ is a constant.
From (\ref{3.5}) it follows  that
$q=\displaystyle{
\frac{
\overline{a}-
\underline{a}}
{\overline{a}+
\underline{a}}}$,
where
$\underline{a}:=\min\limits_{x\in \Omega}{a}(x)$
and
$\overline{a}:=\max\limits_{x\in \Omega}{a}(x)$. Then $\rho_*=\frac{2 a_0}{\overline {a}+\underline{a}}$ and the iteration procedure (\ref{2.9})
with $\rho=\rho_*$ has the form
\ben
\label{3.6}
\IntO Qu_k\odot Q w\,dx=
%\IntO Q u_{k-1}\odot Q w\,dx-
%\frac{2}{\overline {a}+\underline{a}}\left(\IntO a(x) Q u_{k-1}\odot Q w\,dx-\IntO fw\,dx
%\right)=\\
\IntO\left(1-\frac{2a}{\overline {a}+\underline{a}}\right) 
Q u_{k-1}\odot Q w\,dx+\frac{2}{\overline {a}+\underline{a}}\IntO fw\,dx
\een
 From Theorem \ref{th1}, it follows that
\be
\IntO |Q(u_k-u)|^2\,dx\,\leq  C\left(\displaystyle{
\frac{
\overline{a}-
\underline{a}}
{\overline{a}+
\underline{a}}}\right)^{2k}.
\ee
{\it Example 2. Oscillation around a given function.}
Consider a somewhat different example.
Let $a(x)$ be a function
oscillating around a certain mean function $g(x)$ so that 
\be 
\displaystyle\frac{a(x)}{g(x)}\in [1-\epsilon,1+\epsilon],
\quad 
\epsilon\in (0,1).
\ee
 If $g$ is a relatively simple
function, then it is natural to set
$a_\circ(x)=g(x)$. 
By (\ref{3.4}), we find that
$
{\mathds h}_\oplus=\,
1+\epsilon$,
$
{\mathds h}_\ominus=1-\epsilon
$, and $q=\epsilon$.
Hence the method is very efficient for small $\epsilon$
(i.e., if $a$
oscillates around $g$ with a relatively small amplitude).
Figures \ref{fig:1DPeriodStruct2} and \ref{fig:2DPeriodStruct2} illustrate three examples of 
quasi-periodic coefficients  
$a$ and respective $a_\circ$ corresponding to the case of 
oscillation around constant with smooth modulation, 
oscillation around given smooth function, 
or oscillation around piecewise constant function.

{\it Example 3. Piecewise constant coefficients.}
Consider a more complicated case, where $\Omega$ is divided into
$N$ nonoverlapping subdomains $\Omega_i$ and $\Lambda_\circ(x)=c_i{\mathbb I}$ if
 $x\in \Omega_i$. 
Define the numbers $a^{(i)}_\oplus:=\max\limits_{x\in \Omega_i} a(x)$,
 $a^{(i)}_\ominus:=\min\limits_{x\in \Omega_i} a(x)$,
\be 
{\mathds h}_\ominus=\min\left\{\frac{a^{(1)}_\ominus}{c_1},\frac{a^{(2)}_\ominus}{c_2},...,\frac{a^{(N)}_\ominus}{c_N}\right\},\quad{\rm and}\quad
{\mathds h}_\oplus=\max\left\{\frac{a^{(1)}_\oplus}{c_1},\frac{a^{(2)}_\oplus}{c_2},...,\frac{a^{(N)}_\oplus}{c_N}\right\}.
\ee
Since the constantans $c_i$ are defined up to a common multiplier,
we can without a loss of generality assume that
\ben
\label{3.7}
\sum^{(N)}_{i=1}\lambda_i=1, \qquad{\rm where}\qquad \lambda_i=\frac{1}{c_i}.
\een
In accordance with (\ref{3.5}),
maximum of ${\mathcal Q}(\Lambda,\Lambda_\circ)$ is attained
if
\ben
\label{3.8}
%{\mathcal F}(\lambda_1,\lambda_2,...,\lambda_N):=
\frac{\min\left\{\lambda_1 a^{(1)}_\ominus,\lambda_2 a^{(2)}_\ominus,...,\lambda_N a^{(N)}_\ominus\right\}}
{\max\left\{\lambda_1 a^{(1)}_\oplus,\lambda_2 a^{(2)}_\oplus,...,\lambda_N a^{(N)}_\oplus\right\}}\,\rightarrow\,max,
\een
where $\lambda_i>0$ and satisfy (\ref{3.7}).
If $N=2$, then the problem (\ref{3.8}) has a simple solution, which shows that the ratio
$\frac{\lambda_1}{\lambda_2}$ (i.e., $\frac{c_2}{c_1}$) can be any in the interval
$[\xi_1,\xi_2]$, where
$\xi_1=\min\{\frac{a^{(2)}_\ominus}{a^{(1)}_\ominus},  \frac{a^{(2)}_\oplus}{a^{(1)}_\oplus} \} $ and
$\xi_2=\max\{\frac{a^{(2)}_\ominus}{a^{(1)}_\ominus},  \frac{a^{(2)}_\oplus}{a^{(1)}_\oplus} \} $.

It is interesting to compare these results with those generated
by homogenized models in the case of perfectly periodic structures. For this purpose, we consider a simple $1$-dimensional problem
\be 
(a u^\prime)^\prime-f=0\qquad{\rm in}\;(0,1)
\ee
with
\be 
&&a(x)=a^{(1)}(x)\qquad {\rm in}\;\Omega_1=(0,\beta),\qquad\beta\in (0,1),\\
&&a(x)=a^{(2)}(x)\qquad {\rm in}\;\Omega_2=(\beta,1),
\ee
where $a^{(1)}(x)$ is a perfectly periodical function
attaining only two values $a^{(1)}_\oplus$ (Lebesgue measure
of this set is $\kappa_1|\Omega_1|$, $\kappa_1\in (0,1)$) and 
$a^{(1)}_\ominus$ (Lebesgue measure
of this set is $(1-\kappa_1)|\Omega_1|$). Similarly, 
$a^{(2)}(x)$ is a perfectly periodical function
attaining only two values $a^{(2)}_\oplus$ (Lebesgue measure
of this set is $\kappa_2|\Omega_2|$, $\kappa_2\in (0,1)$) and 
$a^{(2)}_\ominus$ (Lebesgue measure
of this set is $(1-\kappa_2)|\Omega_2|$). Assume that the amount of periods
is very large and, therefore, the homogenization method can be successfully applied. The corresponding  homogenized problem has the following coefficients
\be 
\wh a^{(1)}:=\left(\frac{1}{\beta}\int\limits^{\beta}_0\frac{1}{a^{(1)}(x)}\,dx\right)^{-1}\quad{\rm in}\;\Omega_1\quad{\rm and}\quad
\wh a^{(2)}:=\left(\frac{1}{1-\beta}\int\limits^{1}_{\beta}\frac{1}{a^{(2)}(x)}\,dx\right)^{-1}\quad{\rm in}\;\Omega_2.
\ee
It is easy to see that
\be 
\wh a^{(1)}=\frac{a^{(1)}_\ominus a^{(1)}_\oplus}{\kappa_1 a^{(1)}_\ominus+(1-\kappa_1)a^{(1)}_\oplus}\in (a^{(1)}_\ominus,a^{(1)}_\oplus),\qquad
\wh a^{(2)}=\frac{a^{(2)}_\ominus a^{(2)}_\oplus}{\kappa_2 a^{(2)}_\ominus+(1-\kappa_2)a^{(2)}_\oplus}
\in (a^{(2)}_\ominus,a^{(2)}_\oplus)
\ee
Hence
\be 
\frac{\wh a^{(2})}{\wh a^{(1)}}\in (\zeta_1,\zeta_2), 
\quad{\rm where}\quad
\zeta_1:=\frac{\min\{a^{(1)}_\ominus,a^{(2)}_\ominus\}}
{\max\{a^{(1)}_\ominus,a^{(2)}_\ominus\}},\quad
\zeta_2:=\frac{\max\{a^{(1)}_\ominus,a^{(2)}_\ominus\}}
{\min\{a^{(1)}_\ominus,a^{(2)}_\ominus\}}.
\ee
It is clear that $\zeta_1\leq\xi_1$ and $\zeta_2\geq\xi_2$. 
Therefore, homogenized coefficients {\em may not generate the best
piece wise constant $a_\circ$}, which produces the smallest
contraction factor $q$.

%%%%%%%%%%%%%%%%%%%%%%%%%%%%%%%%%%%%%%%%%

\section{Error estimates}
\subsection{General estimate}
Since $T_\rho$ is a contractive mapping, we can use the Ostrowski estimates (see \cite{Ostrowski,Zeidler,ReGruyter}), which yield the estimate
%\begin{multline*}
%\frac{1}{1+q}\|T_\rho v-v\|_\circ\,\leq\,\|v-u\|_\circ\leq\|v-T_\rho v\|_\circ+
%\| T_\rho v-u\|_\circ\leq\\
%\qquad \leq (1+\frac{q}{1-q})\|v-T_\rho v\|_\circ=
%\frac{1}{1-q}\| T_\rho v-v\|_\circ.
%\end{multline*}
of the distance between $v\in V$ and the fixed point:
\ben
\label{4.1}
\|v-u\|_\circ\,\in
\left\{\frac{\epsilon}{1+q(\rho)},\,\frac{\epsilon}{1-q(\rho)}
\right\},\qquad{\rm where}\;\epsilon:=\|T_\rho v-v\|_\circ.
\een
The is estimate cannot be directly applied because $v_\rho:=T_\rho v$ is
 generally unknown (it is the exact solution of a boundary value problem). 
 Instead, we must use a numerical approximation
$\wt v_\rho$ (in our analysis, we impose no restrictions
on the method by which the function $\wt v_\rho\in V$ was constructed).
Thus, the difference $\eta_\rho:=v-\wt v_\rho$ is a known function
and the quantity $\delta_\rho=\|\eta_\rho\|_\circ$ is directly computable. It is easy to see that
\ben
\label{4.2}
\delta_\rho -\|\wt v_\rho-v_\rho\|_\circ\,\leq\,\|v_\rho-v\|_\circ\,\leq\, \delta_\rho+\|\wt v_\rho-v_\rho\|_\circ.
\een
To deduce a fully computable
majorant of the norm $\|\wt v_\rho-v_\rho\|_\circ$ we use the method suggested
in \cite{Re2000,ReGruyter}.
First, we rewrite (\ref{2.1}) in the form
\ben
\label{4.3}
(\Lambda_\circ Q v_\rho, Qw)=(\Lambda_\circ Qv, Qw) - \rho \Bigl( (\Lambda Q v, Qw)-
<f,w>\Bigr) .%\qquad \forall w\in V
\een
For any $y\in Y$ and $w\in V_0$, we have
\ben
\label{4.4}
&(\Lambda_\circ Q (v_\rho-\wt v_\rho), Qw)&=(\Lambda_\circ Q
( v -\wt v_\rho), Qw))-\rho\,\Bigl( (\Lambda Q v , Qw)-
<f,w>\Bigr)\\
&&=(\Lambda_\circ Q
( v -\wt v_\rho)-\rho\Lambda Q v +y, Qw))-<Q^*y+\rho f,w>. 
\nonumber
\een
We estimate the first term in the right hand side
of (\ref{4.4}) as follows:
\begin{multline*}
(\Lambda_\circ Q
( v -\wt v_\rho)-\rho\Lambda Q v +y, Q( v_\rho -\wt v_\rho)))=
(Q(v -\wt v_\rho)-\rho\Lambda^{-1}_\circ\Lambda Q v +\Lambda^{-1}_\circ y,\Lambda_\circ Q( v_\rho -\wt v_\rho)))\\
\leq 
\left(
\Lambda_\circ Q(v -\wt v_\rho)+\tau, 
Q(v -\wt v_\rho)+\Lambda^{-1}_\circ\tau\right)^{1/2}\|v_\rho -\wt v_\rho\|_\circ,
\end{multline*}
where $\tau:=y-\rho \Lambda Q v $. The second term meets the estimate
\be 
<Q^*y+f,v_\rho-\wt v_\rho>\,\leq\,\NNN Q^*y+\rho f \NNN\,\|v_\rho-\wt v_\rho\|\,
\leq\,\frac{1}{(\laomin)^{1/2}}   \NNN Q^*y+\rho f\NNN\,\| v_\rho -\wt v_\rho\|_\circ,
\ee
where
$
\NNN w^*\NNN\,=\,\sup\limits_{w\in V}\frac{<w^*,w>}{\|w\|_\circ}
$
is the dual norm. Hence,
\ben
\label{4.5}
\| v_\rho -\wt v_\rho\|_\circ\,\leq\,\left(
\Lambda_\circ Q\eta_\rho+\tau, 
Q\eta_\rho+\Lambda^{-1}_\circ\tau\right)^{1/2}+\frac{1}{\sqrt{\laomin}}   \NNN Q^*y+\rho f\NNN=:M_\oplus(\eta_\rho,\tau).
\een
Notice that
\be
\inf\limits_{y\in Y}
M_\oplus(\eta_\rho,\tau)=\|v_\rho-\wt v_\rho\|_\circ.
\ee
Indeed, set $y=\Lambda_\circ Q(v_\rho- v)+\rho \Lambda Q v$.
Then, $\tau=\Lambda_\circ Q(v_\rho-v)$. In view of (\ref{4.3}),
$Q^* y+\rho f=0$,
and the majorant is equal to $\|v_\rho-\wt v_\rho\|^2_\circ$.
Hence, the estimate (\ref{4.5}) has no gap. 

It is worth noting that computation of the majorant $M_\oplus$
does not require inversion of the operator $\Lambda$ associated with a complicated quasi--periodic problem.

\begin{remark} 
 $M_\oplus(\eta_\rho,\tau)$ is an a posteriori error majorant
 of the functional type (its derivation is performed by purely functional
 methods based on generalized formulation of the boundary value problem and
 special properties of approximations or numerical method are not used).
Properties of such type error majorants are well studied
(see \cite{Re2000,ReGruyter}  and the literature cited therein). It is not difficult to
show that the last term of $M_\oplus(\eta_\rho,\tau)$ can be estimated
via an explicitly computable quantity provided that $y$ has the same regularity
as the true flux. However, in our subsequent analysis these advanced forms
of the majorant are not required. Therefore we omit this discussion (interested reader can find
the respective analysis in \cite{ReGruyter}). Numerous
tests performed for different boundary value problems have
confirmed high practical efficiency of error majorants of the functional type.  It was shown that  $M_\oplus$ is a guaranteed and efficient majorant of the global 
error and generates good indicators of local 
errors if $y$ is replaced by a certain numerical reconstruction of the exact
dual solution. 
There are many different ways to obtain suitable reconstructions (see \cite{MaNeRe} for a systematic discussion of
computational aspects of this error estimation
method). 
Error majorants of this type  can be also used for the evaluation of modeling errors (see \cite{ReSaSm,ReSaSa2}).
\end{remark}

Now, (\ref{4.1}), (\ref{4.2}), and (\ref{4.5}) yield the following result
\begin{theorem}
\label{th2}
The error  $e=v-u$ is subject to the
estimate
\ben
\label{4.6}
\|e\|_\circ\,\in
\left[\;\max\left\{0,\,\frac{\delta_\rho-M_\oplus(\eta_\rho,\tau)}{1+q(\rho)}\right\},\;\frac{\delta_\rho+
M_\oplus(\eta_\rho,\tau)}{1-q(\rho)}\;
\right],
\een
where  $\tau:=y-\rho \Lambda Q v $ and $y$ is a function in $Y$
and $M_\oplus$ is defied by (\ref{4.5}).

If $Q^*y+\rho f=0$ then
$
M^2_\oplus(\eta_\rho,\tau)=
 (\Lambda_\circ Q \eta_\rho, Q \eta_\rho)+
(\Lambda_\circ^{-1}\tau, \tau)-2(Q\eta_\rho,\tau ).
$
\end{theorem}

\subsection{Examples}
Now we shortly discuss applications of  Theorem \ref{th2}
to problems, where $Q$ and $Q^*$ are defined by the operators
$\nabla$ and $\dvg$, respectively, $\Lambda_\circ=a_\circ(x){\mathbb I}$,
$\Lambda=a(x){\mathbb I}$, $x\in \Omega$, and $V=\oH^1(\Omega)$.

%If
%$
%\Omega\subset\left\{x\in {\mathbb R}^d\,\mid\,a_s<x<b_s,\;
%b_s-a_s=l_s,\;s=1,2,..,d\right\},
%$
%then
%$C_\Omega=\frac{1}{\kappa\pi}$, $
%\kappa^2=\sum\limits^d_{s=1}\frac{1}{l^2_s}.
%$
%In particular,
%$C_\Omega=\frac{1}{d^{1/2}\pi}$  for $\Omega=(0,1)^d$.

%%%%%%%%%%%%%%%%%%%%%%%%%%%%

%Comment:\\
%\be
%\frac{1}{1+q}\|v_1-v_0\|\,\leq\,\|v_0-v_*\|\leq\|v_0-v_1\|+
%\|v_1-v_*\|\leq\\
%\qquad \leq (1+\frac{q}{1-q})\|v_0-v_1\|=\frac{1}{1-q}\|v_1-v_0\|.
%\ee
%We set $v_0=v$ and $v_1=u$.

\subsubsection{$d=1$}
Let $\Omega=(0,1)$. The equation (\ref{1.1}) 
has the form $(a(x) u^\prime)^\prime-f=0$. 
In this case,  $Qw=w^\prime$,  $Q^*y=-y^\prime$, and (\ref{4.3}) is reduced to
\ben
\label{4.7}
\int\limits^1_0 a_\circ (v_\rho-v)^\prime w^\prime\,dx+\rho\int\limits^1_0(a v^\prime w^\prime+fw)\,dx=0.
\een

In order to apply Theorem \ref{th2},
we set $y=\rho(g(x)+\mu)$, where $g(x)=-\int^x_0fdx$ and $\mu$ is a constant. Then $-y^\prime-\rho f=0$
and $\tau=\rho(g(x)+\mu)-\rho a v^\prime=\rho(\mu+g-a v^\prime)$. The
best constant $\mu$ is defined by minimization of $M^2_\oplus(\eta_\rho,\tau)$, which has the form 
\be
 \int\limits^1_0 (a_\circ(\eta^\prime_\rho)^2+
a_\circ^{-1}\rho^2(\mu+g-a v^\prime)^2-2\eta^\prime_\rho\rho(\mu+g-a v^\prime) dx
\ee
 Since $\int^1_0 \eta^\prime_\rho dx=0$,
%
%and need to find $c$ minimizing the quantity
%\be
%\int\limits^1_0 (a_0|\eta^\prime_{k+1}|^2 +a^{-1}_\circ\rho^2(c-g(x)- a_\epsilon u^\prime_{k,h})^2-2\rho\eta^\prime_{k+1} (c-g(x)- a_\epsilon u^\prime_{k,h})dx.
%\ee
the problem is reduced to
 minimization of the second term
%\be
%\int\limits^1_0 a^{-1}_\circ\rho^2(c-g(x)- a_\epsilon u^\prime_{k,h})^2dx
%\ee
and the best
 $\mu$ satisfies the equation
$
\int\limits^1_0 a^{-1}_\circ
(\mu+g(x)-a v^\prime)dx=0.
$
Hence
\be
\mu=\bar\mu:=\frac{\int^1_0 a^{-1}_\circ(a v^\prime-g)dx}{
\int^1_0 a^{-1}_\circ\,dx},
\ee
and (\ref{4.6}) yields the estimate
\ben
\label{4.8}
\|e\|_\circ\,\in
\left[\;\max\left\{0,\,\frac{\delta_\rho-I_\oplus(v,\wt v_\rho)}{1+q(\rho)}\right\},\;\frac{\delta_\rho+
I_\oplus(v, \wt v_\rho)}{1-q(\rho)}
\right],
\een
where
\be
I^2_\oplus(v, \wt v_\rho)=
\int\limits^1_0 a^{-1}_\circ\Bigl(a_\circ(v-\wt v_\rho)^\prime-
\rho(\bar\mu+g-a v^\prime)\Bigr)^2 dx.
\ee
Here $v$ and $\wt v_\rho$ are two consequent numerical approximations (e.g., finite
element approximations $v^k_h$ and $v^{k+1}_h$ computed
on a mesh ${\mathcal I}_h$.
Then 
\be
\eta_\rho=\eta^k_h:=v^k_h-v^{k+1}_h\;{\rm and}\; \delta_\rho=\delta^k:=\|v^k_h-v^{k+1}_h\|_\circ
\ee
 are directly computable.
 Since $a_\circ$ is a "simple" function, the integrals 
\be
 F_1=\int^1_0\!\!\! a^{-1}_\circ\,dx,\; F_2= \int^1_0\!\!\! a^{-1}_\circ\, g\,dx,\; F_3= \int^1_0\!\!\! a_\circ\left(\eta^{k\,\prime}_h\right)^2\,dx,\;
 F_4= \int^1_0\!\!\! a_\circ\left(\bar\mu+g\right)^2\,dx,\;
 F_5=\int^1_0\!\!\!f\eta^k_h\,dx
\ee
 are easy to compute. Other integrals
\be
 G_1=\int^1_0\!\!\! a^{-1}_\circ a\, v^{k\,\prime}_h\,dx,\;
 G_2=\int^1_0\!\!\! a \,(v^k_h)^\prime\,\eta^{k\,\prime}_h \,dx,\;
 G_3=\int^1_0\!\!\!(\bar\mu+g) a^{-1}_\circ a \,v^{k\,\prime}_h\,dx,\;
 G_4=\int^1_0\!\!\! a^{-1}_\circ \,a^2\,\left(v^{k\,\prime}_h)\right)^2\,dx
\ee
contain highly oscillating coefficient $a$ multiplied by  piece wise polynomial mesh functions.
If $a$ has a low QTT rank tensor representation \cite{KhQuant:09}, then the  integrals can be 
efficiently computed by tensor type methods already discussed in \cite{BokhSRep:15}. 
We have
\be
&&I^2_\oplus(v,\wt v_\rho)= F_3+2 G_2+2\rho F_5
+\rho^2( F_4-2 G_3 + G_4)=:\varepsilon^k,\qquad
\bar \mu=\frac{ G_1- F_2}{ F_1}.
\ee
Here
$\rho=\frac{2}{{\mathds h}_\ominus+{\mathds h}_\oplus}$ is selected in accordance 
with  Section 3. The
respective contraction factor is $q=\frac{{\mathds h}_\oplus-{\mathds h}_\ominus}
 { {\mathds h}_\ominus+{\mathds h}_\oplus}$.
Now (\ref{4.8}) yields easily computable lower and upper bounds of the error
encompassed in $v^k_h$:
\be
\label{4.9}
\frac{\delta^k-\varepsilon^k}{1+q}\,\leq\,\|v^k_h-u\|_\circ\,\leq\, \frac{\delta^k+\varepsilon^k}{1-q}
\ee

%%%%%%%%%%%%%%%%%%%%%%%%%%

\subsubsection{$d=2$}
Computation of $M_\oplus$ for 2d problems can be also
reduced to the computation of one dimensional integrals.
Certainly on the multidimensional case the amount of 
integrals is much larger. However the basic tensor decomposition methods remain the same. 
Below we briefly discuss them with the paradigm of a simple case where
 \be
 f=\fx(x_1)\fy(x_2)\quad{\rm  and} \quad
 a=\ax(x_1)\ay(x_2).
 \ee
  Assume that approximations
are represented in the form of series formed by one dimensional
functions $\phix_{i}$ and $\phiy_{j}$ (which may be supported locally or globally), so that
\be
v=\sum\limits^{n_1}_{i=1}\sum\limits^{n_2}_{j=1} \gamma_{i j}\phix_{i}(x_1)\phiy_{j}(x_2),\quad
\wt v_\rho=\sum\limits^{n_1}_{i=1}\sum\limits^{n_2}_{j=1} \wt\gamma_{i j}\phix_{i}(x_1)\phiy_{j}(x_2).
\ee
In this case,
\be
\nabla \eta_\rho=\left( \sum\limits^{n_1}_{i=1}\sum\limits^{n_2}_{j=1}
 \varsigma_{ij}\frac{\partial\phix_{i}}{\partial x_1}\phiy_{j}  \, ,\,
 \sum\limits^{n_1}_{i=1}\sum\limits^{n_2}_{j=1}\varsigma_{ij} 
 \phix_{i}\frac{\partial\phiy_{j}}{\partial x_2}   \right),
 \quad{\rm where}\quad \varsigma_{ij}=\gamma_{ij}-\wt \gamma_{ij}.
\ee
We define another set of one dimensional functions
$\Wx_{k}(x_1)$ and $\Wy_{l}(x_2)$, which form the vector function
\ben
\label{4.9}
y=\Upsilon_0+\sum\limits^{m_1}_{k=1}
\sum\limits^{m_2}_{l=1}\sigma_{kl}\Upsilon_{kl},\quad
\Upsilon_{kl}=\left\{ \Wx_{k} \frac{\partial\Wy_{l}}{\partial x_2} \; ;\;-\frac{\partial\Wx_{k}}{\partial x_1}\Wy_{l}   \right\}.
\een
Here $\Upsilon_0$ is a given function,
which can be defined in different ways. In particular,
we set $\Upsilon_0=\left\{\Wx_0(x_1)\Wy_0(x_2)\,;\,0\right\}$,
$\Wx_0(x_1)=\int^{x_1}_0\fx dx_1$
and $\Wy_0=-\rho\fy$. The functions $\Upsilon_{kl}$
must satisfy the usual linear independence conditions in order to guarantee
unique solvability of the respective approximation problem.
%Notice that
%\be
%\IntO (\Upsilon_0\cdot \nabla w-\rho fw) dx_1dx_2=0.
%\IntO(\Wx_0\Wy_0\frac{\partial w}{\partial x_1}-\rho fw)dx_1dx_2=0.
%\ee
For any smooth function $w$ vanishing on $\partial\Omega$, we
have
\be
\IntO (\Upsilon_0\cdot \nabla w-\rho fw) dx_1dx_2=0\quad
{\rm and}\quad
\IntO \Upsilon_{kl}\cdot \nabla w dx_1dx_2=0.
%\IntO \left(\Wx_{k} (x_1)\frac{\partial\Wy_{l}(x_2)}{\partial x_2}\frac{\partial w}{\partial x_1}-\frac{\partial\Wx_{k}(x_1)}{\partial x_1}\Wy_{l}(x_2)\frac{\partial w}{\partial x_2}\right)dx_1dx_2=0
\ee
Thus,
$\NNN Q^*y+\rho f\NNN=\NNN\dvg y-\rho f\NNN=0$ and we can use the simplified
form of $M_\oplus$.

In the simplest case
$\Lambda_\circ=a_\circ{\mathbb I}$, where $a_\circ$ is a constant. The best $y$  minimizes the quantity
% $\tau:=y-\rho \Lambda Q v $\qquad
%$
%(\Lambda_\circ^{-1}(y-\rho \Lambda Q v), y-\rho \Lambda Q v)-2(Q\eta_\rho, y-\rho \Lambda Q v ).
%$
\begin{multline}
\label{4.10}
M^2_\oplus(\eta_\rho,\tau)=
 \IntO a_\circ \nabla \eta_\rho\cdot \nabla \eta_\rho dx+
 \IntO a^{-1}_\circ y\cdot y dx+\rho^2\IntO a^{-1}_\circ a^2
 \nabla v\cdot \nabla v dx\\
- 2\IntO (\rho a^{-1}_\circ a\nabla v+\nabla \eta_\rho)\cdot y\,dx+
 2\rho\IntO a\nabla \eta_\rho\cdot\nabla \eta_\rho\,dx,
 \end{multline}
 which shows that  $y$ must satisfy the relation $y=\rho a\nabla v+a_\circ\nabla\eta_\rho$.
 We select $\sigma_{kl}$ that defines Galerkin approximation of this function
% \be
% \IntO (y-\rho a\nabla v+a_\circ\nabla\eta_\rho)\cdot
% \Upsilon_{st}\,dx\,=0\qquad \forall \Upsilon_{st},\,
% s=1,2,...,m_1,\,l=1,2,...,m_2
% \ee
  and arrive at the system
 \begin{multline}
 \label{4.11}
 \sum\limits^{m_1}_{k=1}
\sum\limits^{m_2}_{l=1}\sigma_{kl}\IntO \Upsilon_{kl}\cdot\Upsilon_{st}
dx_1dx_2+\IntO \Upsilon_0\cdot\Upsilon_{st} dx_1dx_2\\
=\sum\limits^{n_1}_{i=1}
\sum\limits^{n_2}_{j=1}
\IntO (\rho a \gamma_{ij}+a_\circ \varsigma_{ij})\left(\frac{\partial\phix_{i}}{\partial x_1}\phiy_{j},
 \phix_{i}\frac{\partial\phiy_{j}}{\partial x_2}\right)\cdot \Upsilon_{st}
dx_1dx_2
 \end{multline}
 
 Introduce the following matrixes 
 \be
&{\rm D}^{(1)}=\left\{D^{(1)}_{kl}\right\}
,\quad D^{(1)}_{kl}=\displaystyle\int^a_0\frac{\partial\Wx_k}{\partial x_1}
\frac{\partial\Wx_l}{\partial x_1}dx_1,\quad &{\rm W}^{(1)}=
\left\{W^{(1)}_{kl}\right\},\;
W^{(1)}_{kl}=\int^a_0 \Wx_k\Wx_l\,dx_1,\\
&{\rm D}^{(2)}=\left\{D^{(2)}_{kl}\right\}
,\quad D^{(2)}_{kl}=\displaystyle\int^b_0\frac{\partial\Wy_k}{\partial x_2}
\frac{\partial\Wy_l}{\partial x_2}dx_2,\quad &{\rm W}^{(2)}=
\left\{W^{(2)}_{kl}\right\},\quad
W^{(2)}_{kl}=\int^b_0 \Wy_k\Wy_l\,dx_2,
\ee
\be
&{\rm F}^{(1)}=\left\{F^{(1)}_{ik}\right\}
,\quad F^{(1)}_{ik}=\displaystyle\int^a_0\frac{\partial\phix_i}{\partial x_1}
\Wx_k dx_1,\quad &{\rm G}^{(1)}=
\left\{G^{(1)}_{ik}\right\},\quad
G^{(1)}_{ik}=
\int\limits^a_0 \phix_i\frac{\partial\Wx_k}{\partial x_1}\,dx_1,\\
&{\rm F}^{(2)}=\left\{F^{(2)}_{jl}\right\}
,\quad F^{(2)}_{jl}=\displaystyle\int^b_0\phiy_j\frac{\partial\Wy_l}{\partial x_2}
 dx_2,\quad &{\rm G}^{(2)}=
\left\{G^{(2)}_{jl}\right\},\quad
G^{(2)}_{jl}=
\int\limits^b_0 \frac{\partial\phiy_j}{\partial x_2}\Wx_l\,dx_2,
\ee
\be
&\wh{\rm F}^{(1)}=\left\{\wh F^{(1)}_{ik}\right\},\; \wh F^{(1)}_{ik}=\displaystyle\int^a_0a_1(x_1)\frac{\partial\phix_i}{\partial x_1}
\Wx_k dx_1,\; &\wh{\rm G}^{(1)}=
\left\{\wh G^{(1)}_{ik}\right\},\;
\wh G^{(1)}_{ik}=
\int^a_0 a_1(x_1) \phix_i\frac{\partial\Wx_k}{\partial x_1}\,dx_1,\\
&\wh{\rm F}^{(2)}=\left\{\wh F^{(2)}_{jl}\right\},\; \wh F^{(2)}_{jl}=\displaystyle\int^b_0a_2(x_2)\phiy_j\frac{\partial\Wy_l}{\partial x_2}
 dx_2,\; &\wh{\rm G}^{(2)}=
\left\{\wh G^{(2)}_{jl}\right\},\;
\wh G^{(2)}_{jl}=
\int^b_0a_2(x_2) \frac{\partial\phiy_j}{\partial x_2}\Wx_l\,dx_2.
\ee
and vectors
\be
\bg^{(1)}=\{g^{(1)}_k\},\quad g^{(1)}_k=\displaystyle\int^a_0 \Wx_0\Wx_k\,dx_1,\qquad
\bg^{(2)}=\{g^{(2)}_l\},\quad g^{(2)}_l=
\displaystyle\int^b_0 \Wy_0\frac{\partial\Wy_l}{\partial x_2}\,dx_2.
\ee
Notice that all  coefficients are presented
by one dimensional integrals, which can be efficiently computed
with the help of special (tensor type) methods (see, e.g.,
\cite{KhQuant:09}-\cite{KhSautVeit:11}).

It is not difficult to see that
\be
Y_{klst}:=\IntO \Upsilon_{kl}\cdot\Upsilon_{st}\,dx
%=\IntO \left\{ \frac{\partial\Wy_{l} }{\partial x_2}\Wx_{k}   \; ;\;-\Wy_{l} \frac{\partial\Wx_{k} }{\partial x_1}   \right\}\cdot
%\left\{ \frac{\partial\Wy_{t} }{\partial x_2}\Wx_{s}   \; ;\;-\Wy_{t} \frac{\partial\Wx_{s} }{\partial x_1}   \right\}=\\
%\IntO \Bigl(\Wx_{k}  \Wx_{s}  
%\frac{\partial\Wy_{l} }{\partial x_2}
%\frac{\partial\Wy_{t} }{\partial x_2} +\Wy_{l} \Wy_{t} \frac{\partial\Wx_{k} }{\partial x_1}   \;\frac{\partial\Wx_{s} }{\partial x_1}   \Bigr)dx\\
=W^{(1)}_{ks}D^{(2)}_{lt}+D^{(1)}_{ks}W^{(2)}_{lt}
\ee
and
\be
\IntO \Upsilon_0\cdot\Upsilon_{st} dx_1dx_2=
\IntO \Wx_0 \Wx_{s}  
\Wy_0 \frac{\partial\Wy_{t} }{\partial x_2}dx_1dx_2=g^{(1)}_{s}
g^{(2)}_{t},
\ee
where ${\rm Y}=\{Y_{klst}\}$ is the fourth order tensor.
Hence the left hand side of the system (\ref{4.11}) has the form ${\rm Y}{\boldsymbol\sigma}+\bg^{(1)}\otimes \bg^{(2)}$.
%\be
%\sum\limits^{m_1}_{k=1}\sum\limits^{m_2}_{l=1}
%Y_{klst}\sigma_{kl}+g^{(1)}_{s}
%g^{(2)}_{t}.
%\ee
In  the right hand side we have the term
\be
\IntO a_\circ \varsigma_{ij}\left(\frac{\partial\phix_{i}}{\partial x_1}\phiy_{j},
 \phix_{i}\frac{\partial\phiy_{j}}{\partial x_2}\right)\cdot \Upsilon_{st}
dx_1dx_2
%=\\
%\IntO a_\circ\varsigma_{ij}\left\{ \frac{\partial\phix_{i}}{\partial x_1}\phiy_{j},
% \phix_{i}\frac{\partial\phiy_{j}}{\partial x_2}  \right\}\cdot
%\left\{ \frac{\partial\Wy_{t} }{\partial x_2}\Wx_{s}   \; ;\;-\Wy_{t} \frac{\partial\Wx_{s} }{\partial x_1}   \right\}=\\
%\IntO a_\circ\varsigma_{ij} \left( \frac{\partial\phix_{i}}{\partial x_1} \Wx_{s}  \phiy_{j} \frac{\partial\Wy_{t} }{\partial x_2} -
% \phix_{i}\frac{\partial\Wx_{s} }{\partial x_1}\,\frac{\partial\phiy_{j}}{\partial x_2} 
% \Wy_{t}    \right)
% =a_\circ(F^{(1)}_{is}
% F^{(2)}_{jt}-G^{(1)}_{is}G^{(2)}_{jt})\varsigma_{ij}
 =a_\circ {\rm H}\bvarsigma,
\ee
where ${\rm H}=\{H_{ijst}\},\,
H_{stij}=
F^{(1)}_{is}F^{(2)}_{jt}-G^{(1)}_{is}G^{(2)}_{jt}$.
Another term is
\be
\IntO \rho a \gamma_{ij}
\left(\frac{\partial\phix_{i}}{\partial x_1}\phiy_{j},
 \phix_{i}\frac{\partial\phiy_{j}}{\partial x_2}\right)\cdot \Upsilon_{st}
dx_1dx_2
%=(\wh F^{(1)}_{is}\wh F^{(2)}_{jt}-\wh G^{(1)}_{is}\wh G^{(2)}_{jt})\gamma_{ij}
= \wh {\rm H}\bgamma,
\ee
where $\wh{\rm H}=\{\wh H_{ijst}\},\,
\wh H_{stij}=
\wh F^{(1)}_{is}\wh F^{(2)}_{jt}-\wh G^{(1)}_{is}\wh G^{(2)}_{jt}$.

Now (\ref{4.11}) implies
${\boldsymbol\sigma}={\rm Y}^{-1}(\wh {\rm H}\bgamma+a_\circ {\rm H}\bvarsigma-\bg^{(1)}\otimes \bg^{(2)})$
%\begin{multline}
%\label{4.12}
%\sum\limits^{m_1}_{k=1}\sum\limits^{m_2}_{l=1}
%\sigma_{kl}(W^{(1)}_{ks}D^{(2)}_{lt}+D^{(1)}_{ks}W^{(2)}_{lt})=\\
%\sum\limits^{n_1}_{i=1}\sum\limits^{n_2}_{j=1}
% \left((\wh F^{(1)}_{is}\wh F^{(2)}_{jt}-\wh G^{(1)}_{is}\wh G^{(2)}_{jt})\gamma_{ij}
%+a_\circ(F^{(1)}_{is}F^{(2)}_{jt}-G^{(1)}_{is}G^{(2)}_{jt})\varsigma_{ij}\right)-g^{(1)}_{s}g^{(2)}_{t}
%\end{multline}
and the value of $M_\oplus$ is obtained by (\ref{4.6}), (\ref{4.9}), and (\ref{4.10}).

%%%%%%%%%%%%%%%%%%%%%%%%%%%%%%%%%%%%%%%%%%%%%%%%%

\color{black}
\section{Low-rank solution of the discrete equation}\label{sec:LowRank_Solut}
%==============\\
% \cred
% Here we should put some text connecting two parts.
% What disturbs me is that in Sect. 5 there are no words about "simplified"
% (preconditioned) operator $A_\circ$. However, this is the main idea exposed
% in Sect. 1-4. May be here it is simply hidden in the algebraic treatment
% of the system $Au=f$. But anyway we should extract all this and explain
% things within the general framework. 
% Below is a short draft where I tried to explain this.
% But there appeared questions. Please could you explain/add/correct!
% \cn

%\color{green}
In what follows 
we assume that $f$ and $a$ admit low rank representation 
(e.g., $f=\sum_{i=1}^{R_f}f^i_1(x_1)f^i_2(x_2)$, $a=\sum_{j=1}^{R_a}a^j_1(x_1)a^j_2(x_2)$).
Then one may assume that the exact FEM solution can be well approximated by
$u^K(x)=\sum^K_{j=1} u^j_1(x_1) u^j_2(x_2)$,
where $K$ depends on the separation rank of $f$ and $a$. In some cases this important property
can be rigorously proven (say, for Laplacian like operators). The similar low rank
approximation can be observed for the QTT tensor approximation (see \cite{BokhSRep:15}). 
Existence of low rank solution means that for some $K$ we have $u_K\approx u$ up to
the rank truncation threshold.
%\color{black}

%\color{blue}
Here we sketch the rank-structured computational scheme.
In our set of examples the original problem: find $u$ such that
\ben
\label{a}
\IntO a(x) \nabla u\cdot \nabla w\,dx=\IntO f w dx\qquad\forall w\in V_0:=H^1_0
\een
is replaced by the Galerkin problem for low rank representations
\ben
\label{b}
\IntO a(x) \nabla u^K\cdot \nabla w^K\,dx=\IntO f w dx\qquad\forall w^K\in V^K_0,
\een
where $V^K_0$ is a subset of $V_0$ formed by functions of the type
\be
w^K(x)=\sum^K_{j=1} \phi^j_1(x_1) \phi^j_2(x_2).
\ee
Therefore, in terms of the general scheme exposed in the introduction, the
Problem $\mathcal P$ is now the problem (\ref{b}) and we solve it by iterations with the help of
simplified (preconditioned) problem
\ben
\label{c}
\IntO a_\circ(x) \nabla u_k\cdot \nabla w^K\,dx=\IntO  f_{k-1} w dx\qquad\forall w^K\in V^K_0,
\een
where $a_\circ$ is a simple (mean) function and $f_{k-1}$ depends on $u_{k-1}$.

% MAIN QUESTION:\\
%  why this discrete preconditioned problem is simpler than the
% discrete problem (\ref{b}) with complex coefficient?
% 
% My guess is may be the answer is as follows:\\
Given the right-hand side, the problem (\ref{c}) is much simpler than the initial equation since
the matrix $\Lambda_\circ$, generated by the coefficient $a_\circ$ is easily invertible.
Moreover, 
the coefficient $a$ may be rather complicated and admits a representation with rank $R$, i.e.,
\be
a(x)=\sum^R_{s=1} a^{(s)}_1(x_1)...a^{(s)}_d(x_d),
\ee
where $R$ is a small integer.
When we construct the low-rank Kronecker representation of stiffness
matrix for this $a$, which is presented by  elements of $4R$ matrices
computed by only 1D integrals containing oscillating functions
$a^{(s)}_i(x_i)$. 

% we must perform calculations as in (5.5) for each term of $a$. 
% Then we need to compute elements of $4N$ matrixes, 
% which are presented by rather complicated integrals containing oscillating functions
% $a^{(s)}_i(x_i)$. 

If we use (\ref{c}), then $a_\circ$
is a simple function, it may be a even a constant, or
a function representable in the form $a^{\circ}_1(x_1)...a^{\circ}_d(x_d)$ 
with very simple multipliers.
Then, the respective Kronecker stiffness matrix $\Lambda_\circ$ is computed
much easier and has a simple (low rank) form that allows the low rank
representation of its inverse.

% In fact the text and formulas (5.5) related to $a_{ij}=...$ and later on 
% should be  addressed to this case.
% May be more logical is to write formulas there in terms of $a_\circ$ instead of $a$
% and show that we obtain a simple Kronecker
% representation of the form $A_1\otimes M_2+M_1\otimes A_2$
% for which we have a very efficient solver based
% on QTT?
%\cn

\subsection{Kronecker product representation of the stiffness matrix}\label{ssec:Kron_matr}

% \color{red} Comment: In new notation this "complicated" coefficient is denoted
% not $a_\epsilon$ but $a$, also $u_\epsilon$ is replaced by $u$ OK, follow on \color{black}\\
% \vspace{0.1cm}

We consider the elliptic diffusion equation  with
quasi-periodic coefficient $a(x)>0$ (whose oscillations are characterized 
by the parameter $\epsilon$) % $x\in [0,1]^2$,
\ben \label{eqn:5.1}
{\mathcal A}u  =-\mbox{div} (a(x) \nabla u) =f(x), 
\quad x=(x_1,\ldots,x_d)\in \Omega=(0,1)^d,\quad {u}_{|\Gamma}=0,
\een
where the function $f$ corresponds to the modified right hand side in the problem (\ref{4.3}),
$\Gamma=\partial \Omega$, and the right-hand side $f(x_1,\ldots,x_d)$ 
can be represented with a low separation rank.

% \color{red} 
% Chto takoe malenkie stolbiki na etih kartinkah? 6 bolshih razeleni na neskolko
% malenkih. Chitatel' udivliaetsa zachem? Esli kazhdii bolshoi stolbik ne sploshnoi a tam toghe nabor razdelennih malenkih stolbikov, to etogo na kartinke ne vidno!\\
% Drugoe zamechanie: stolbiki ne dolzhni opuskatsa do nulia (mi eto obsughdali).
% Inache zadacha budet virozhdennoi!
%  \cn

\begin{figure}[htbp]
\centering
\includegraphics[width=7.0cm]{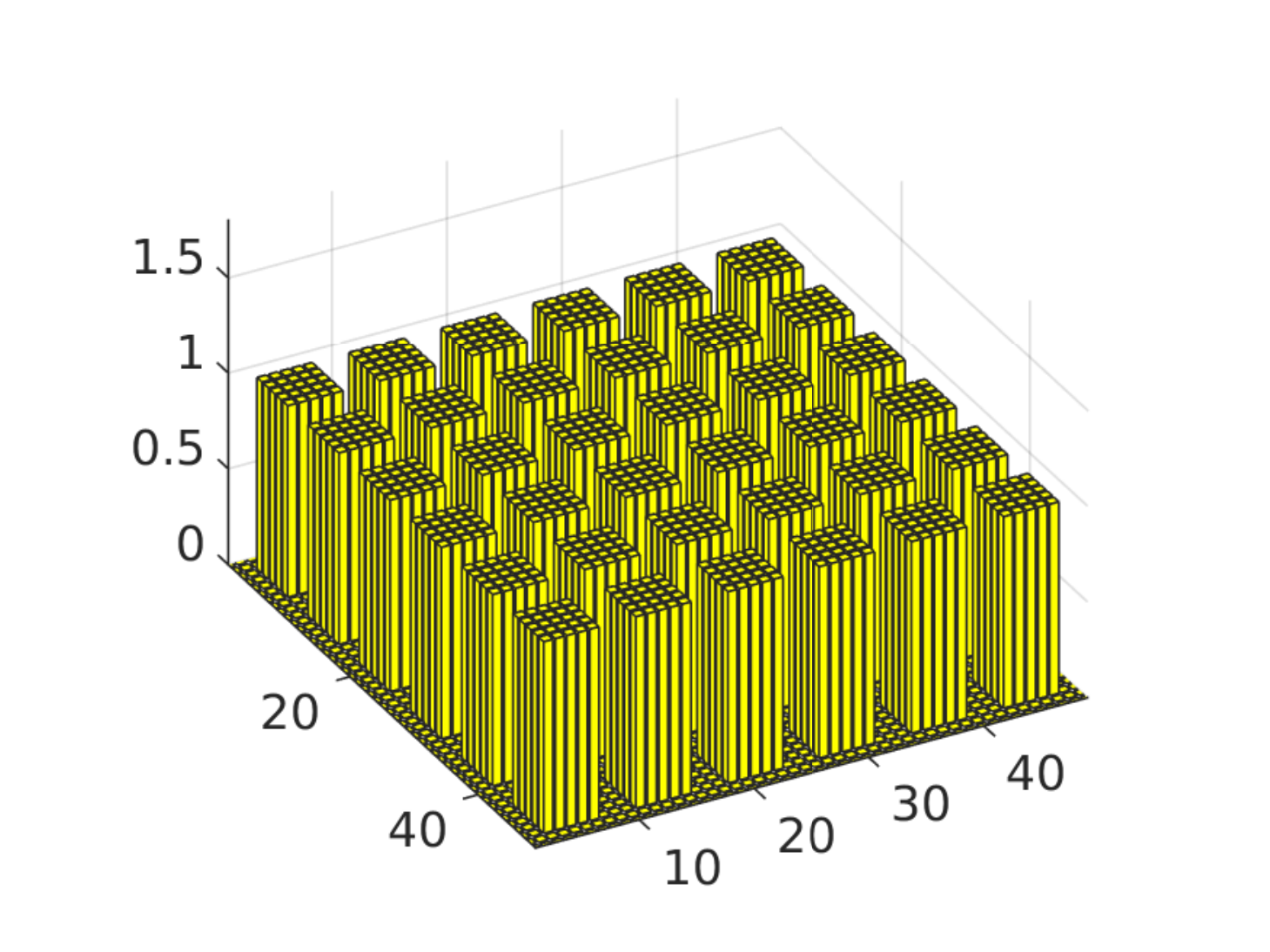}
\includegraphics[width=7.0cm]{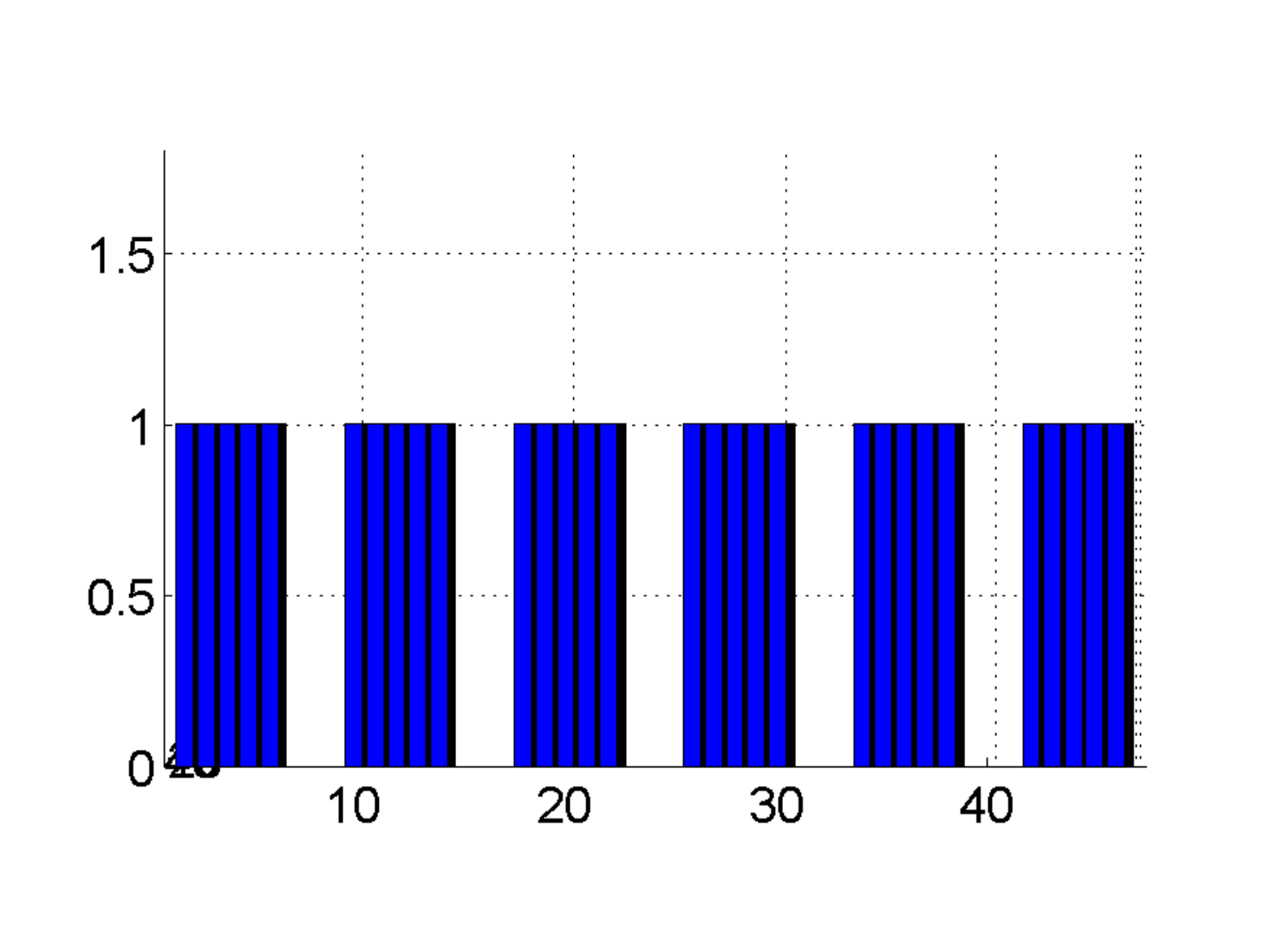}
\caption{\small Example of the 2D periodic oscillating coefficients (left) 
and the 1D factor $a_1(x_1)$.}
\label{fig:Hom2d_Exm1}  
\end{figure}

% \color{green}
Figure \ref{fig:Hom2d_Exm1} illustrates a 2D example of $L \times L$ periodic
coefficient with $L=6$ corresponding to the choice $\epsilon=1/L$. 
In this example, the scalar coefficient 
is represented by the separable function $a(x)=C+ a_1(x_1) a_1(x_2)$, $C>0$,
where the generating univariate function $a_1(x_1)$ has the shape of six uniformly 
distributed bumps of hight $1$ as shown in Figure \ref{fig:Hom2d_Exm1}, right.
Figure \ref{fig:Hom2d_Exm1}, left, presents the oscillating part of 2D coefficients 
function, $a_1(x_1) a_1(x_2)$.

The examples of other possible shapes of the equation coefficient corresponding to
the cases (1), (2) and (3) specified in Introduction are presented in
Figures \ref{fig:1DPeriodStruct2} and \ref{fig:2DPeriodStruct2}.

We apply the FEM Galerkin discretization of equation (\ref{eqn:5.1}) by means of
tensor-product piecewise affine basis functions  (instead of "linear finite elements") 
\be
&&\{\varphi_{\bf i}(x):=\varphi_{i_1}(x_1) \cdots \varphi_{i_d}(x_d)\}, 
\quad {\bf i}=(i_1,\ldots,i_d), \; 
%\\
% &&  \text{\cred are these $\varphi_{i_k}$ are 1D finite element basis functions \cn} \\
%&&
i_\ell\in {\mathcal I}_\ell=\{1,\ldots,n_\ell\}, \; \ell=1,\ldots,d,
\ee
% \cred Let us explicitly write that we have  $N=n_1 n_2...n_d$
% basic functions $\varphi_i$.\cn
%\color{green}
where $\varphi_{i_k}$ are 1D finite element basis functions 
(say, piecewise linear hat functions).
%\cn

We associate the univariate basis functions with the uniform 
grid $\{\nu_j\}$, $j=1,\ldots,n_\ell,$ on $[0,1]$ with the mesh size $h=1/(n_\ell+1)$.
In this construction we have  $N=n_1 n_2...n_d$ basis functions $\varphi_{\bf i}$.
Notice that the univariate grid size $n_\ell$ is of the order of $n_\ell=O(1/\epsilon)$
designating the total problem size $N=O(1/\epsilon^d)$.

% \cred Comment: So far we have not said anything about the grid (probably here
% it concerns a one dimensional grid for the function like $\phi_{i1}$).
% This is an important point.
% I suggest to write about this with more details and later.  DONE \cn
%$i=(i_1,i_2)$ for $i_\ell\in {\mathcal I}_\ell=\{1,\ldots,n_\ell\}$, $\ell=1,\ldots,d$.

For ease of exposition we, first, consider the case $d=2$, and
further assume that the scalar diffusion coefficient $a(x_1,x_2) $ can be 
represented  in the  form 
\[
 a(x_1,x_2) = \sum\limits^{R}_{k=1} a^{(1)}_k (x_1) a^{(2)}_k (x_2) > 0
\]
with a small rank parameter $R$.

The $N\times N$ stiffness matrix is constructed by the standard mapping of the multi-index ${\bf i}$
into the $N$-long univariate index $i$ representing all degrees of freedom. 
For instance,  we use the so-called big-endian convention for $d=3$ and $d=2$
\[
 {\bf i}\mapsto i:= i_3 + (i_2-1)n_3 + (i_1-1)n_2 n_3, \quad {\bf i}\mapsto i:= i_2 + (i_1-1)n_2,
\]
respectively. Hence all matrices and vectors are defined on the long index $i$ as usual,
however, the special Kronecker structure allows the low-storage and low-complexity
matrix vector multiplications when appropriate, i.e. when a vector also admits the 
low-rank Kronecker form representation. 
 In particular, the basis function $\varphi_{\bf i}$ is designated 
via the long index, i.e. $\varphi_i=\varphi_{\bf i}$.

First, we consider the simplest case $R=1$ and let $d=2$.
We construct the Galerkin stiffness matrix $A=[a_{ij}]\in \mathbb{R}^{N\times N}$
in the form of a sum of Kronecker products of small "univariate" matrices.
Recall that given $p_1\times q_1$ matrix $A$ and 
$p_2\times q_2$ matrix $B$, their Kronecker product is defined as a $p_1 p_2 \times q_1 q_2$ 
matrix $C$ via the block representation 
$$ 
C=A \otimes B =[a_{ij}B], \quad i=1,\ldots,p_1, \; j=1,\ldots,q_1.
$$
We say that the Kronecker rank of the matrix $A$ in the representation above equals to $1$.
Now the elements of Galerkin stiffness matrix take a form 
 
% \cred 
% I do not understand the relations below. what is here $i,j$? multi-indexes? 
% what is $\phi_i$? may be $\phi_{\bf i}$ as on the previous page
% where $\bf i$ is the multi-index?
% \cn

\begin{multline}
%\cred a_{\bf ij}?\cn\;
a_{ij}   =  \langle {\mathcal A} \varphi_i, \varphi_j \rangle =
\IntO  a^{(1)} (x_1) a^{(2)} (x_2) 
\nabla \varphi_i (x)\dot \nabla \varphi_j (x) d x \\
    =  {\int\limits^1_0  } a^{(1)} (x_1) 
    \frac{\partial \varphi_{i_1} (x_1)}{\partial x_1 }   
   \frac{\partial \varphi_{j_1} (x_1)}{\partial x_1 } d x_1   
    \int\limits^1_0 a^{(2)} (x_2) \varphi_{i_2}(x_2)\varphi_{j_2}(x_2) d x_2    \\
    +  \int\limits^1_0  a^{(1)} (x_1) \varphi_{i_1}(x_1) \varphi_{j_1} (x_1)d x_1  
   \int\limits^1_0  a^{(2)} (x_2) \frac{\partial \varphi_{i_2} (x_2)}{\partial x_2 } 
   \frac{\partial \varphi_{j_2} (x_2)}{\partial x_2 } d x_2 , 
 \end{multline}
which leads to the rank-2 Kronecker product representation 
\[
 {A} = [a_{ij}]= {A}_1 \otimes M_2 + M_1 \otimes {A}_2,
\]
where $\otimes$ denotes the conventional Kronecker product of matrices.
Here ${A}_1=[a_{i_1j_1}]\in \mathbb{R}^{n_1\times n_1} $ and 
${A}_2=[a_{i_2 j_2}]\in \mathbb{R}^{n_2\times n_2}$ denote the univariate stiffness matrices
and $M_1=[m_{i_1j_1}]\in \mathbb{R}^{n_1\times n_1}$ and 
$M_2=[m_{i_2j_2}]\in \mathbb{R}^{n_2\times n_2}$ define the corresponding 
weighted mass matrices, e.g.,
\[
 a_{i_1j_1}= 
 {\int\limits^1_0  } a^{(1)} (x_1) 
    \frac{\partial \varphi_{i_1} (x_1)}{\partial x_1 }   
   \frac{\partial \varphi_{j_1} (x_1)}{\partial x_1 } d x_1 , \quad
  m_{i_1j_1}= 
  \int\limits^1_0  a^{(1)} (x_1) \varphi_{i_1}(x_1) \varphi_{j_1} (x_1)d x_1 .
\]
By simple algebraic transformations (e.g. by lamping of the tri-diagonal 
mass matrices, which does not effect
the approximation order of the FEM discretization)
the matrix ${A}$ can be simplified to the form   %\cite{..}
\ben \label{eqn:Lapl_Kron_D}
 {A} \mapsto A = A_1 \otimes D_2 + D_1 \otimes A_2,
 \een
 where $D_1, D_2$ are the diagonal matrices.  
 The matrix $A$ corresponds to the FEM discretization of the 
 initial elliptic PDE with complicated highly oscillating coefficients. 
% \cred this I do not understand COMMENT ABOVE.\cn.

%\color{green}
The simple choice of the spectrally equivalent preconditioner $A_\circ$ 
corresponds to the operator Laplacian. In this case 
the representation in (\ref{eqn:Lapl_Kron_D}) is simplified 
 to the discrete Laplacian matrix in the form of rank-$2$ Kronecker sum
\ben \label{eqn:Lapl_Kron}
 {A} \mapsto \Lambda_\circ = A_1 \otimes I_2 + I_1 \otimes A_2,
\een
where $I_1$ and $I_2$ denote the identity matrices of the corresponding size.
This matrix will be used in what follows as a prototype preconditioner for solving
the linear system of equations 
\ben \label{eqn:FEM_syst}
 A {\bf u} = {\bf f}.
\een
The matrix $A$ is constructed in general for the $R$-term separable coefficient 
$a(x_1,x_2)$ with $R\geq1$ which leads to the rank-$2R$ Kronecker sum representation
\[
  A = \sum\limits^{R}_{k=1} [A_{1,k} \otimes D_{2,k} + D_{1,k} \otimes A_{2,k}],
\]
with matrices of the respective size.
%\cn 

% \cred all this above is very unclear. What are the matrixes $I_1$ and $I_2$? 
% How this is related to $A_\circ$?
% Weak place for a reviewer...
% \cn

\subsection{Existence of the low-rank solution}\label{ssec:LowRankSolut}

%\color{green}
In this paper we discuss the approach based on the low rank
separable $\epsilon$-approximation of the solution to the equation (\ref{eqn:FEM_syst})
that is considered as the $d$-dimensional real valued array, 
${\bf u}\in \mathbb{R}^{n_1\times \cdots \times n_d}$. 
In general, for the case $R>1$ this favorable property is not guaranteed by the low Kronecker rank
representation to the Galerkin system matrix $A$, discussed in the previous section.
%\cn 

Let $R=1$ and $d=2$, the existence of the low rank approximation  
to the solution of the equation (\ref{eqn:FEM_syst}) with the low-rank right-hand side 
\[
 {\bf f}=\sum\limits^{R_f}_{k=1} {\bf f}^{(1)}_k \otimes {\bf f}^{(2)}_k,  \quad 
 {\bf f}^{(\ell)}_k\in \mathbb{R}^{n_\ell},
\]
and with the system matrix in the form (\ref{eqn:Lapl_Kron}) 
can be justified by plugging the representation (\ref{eqn:Lapl_Kron}) 
in the $\mbox{sinc}$-quadrature approximation to the Laplace integral transform \cite{GaHaKh4:02} 
\ben \label{eqn:Inv_Kron}
 \Lambda_\circ^{-1}=\int_{\mathbb{R}_+} e^{- t \Lambda_\circ } dt
 \approx B_M:=\sum\limits_{k=-M}^M c_k e^{-t_k \Lambda_\circ }=
 \sum\limits_{k=-M}^M c_k e^{-t_k A_1 } \otimes e^{-t_k A_2 },
\een
taking into account that the matrices $A_1$ and  $A_2$ commute with $I_1$ and $I_2$, respectively.
Hence, the equation (\ref{eqn:Inv_Kron}) represents the accurate rank-$(2M+1)$ Kronecker 
product approximation to $\Lambda_\circ^{-1}$
which can be applied directly to the right-hand side to obtain
\[
 {\bf u} = \Lambda_\circ^{-1}{\bf f} \approx B_M {\bf f} =
 \sum\limits_{k=-M}^M c_k \sum\limits^{R_f}_{m=1} 
 e^{-t_k A_1 }{\bf f}^{(1)}_m \otimes e^{-t_k A_2 }{\bf f}^{(2)}_m.
\]

The numerical efficiency of the representation (\ref{eqn:Inv_Kron}) can be explained by 
the fact that
the quadrature parameters $t_k,c_k$ can be chosen in such a way that the low Kronecker
rank approximation ${B}_M$ converges to $\Lambda_\circ^{-1}$ exponentially fast in $M$.
For example, under the choice $t_k=e^{k h}$, $c_k=h t_k$ with $h=\pi /\sqrt{M}$ there holds
\cite{GaHaKh4:02}
\[
 \|\Lambda_\circ^{-1} - {B}_M \| \leq  C e^{- \beta \sqrt{M}}\|\Lambda_\circ^{-1}\|,
\]
which means that the approximation error $\epsilon>0$ can be achieved with the 
number of terms $R_B=2M+1$ of the order of $R_B=O(|\log \epsilon |^2)$.

\begin{figure}[htbp]
\centering
% %\includegraphics[width=4.0cm]{UexaMUhomo_Sin_G2e14-eps-converted-to.pdf}
\includegraphics[width=5.4cm]{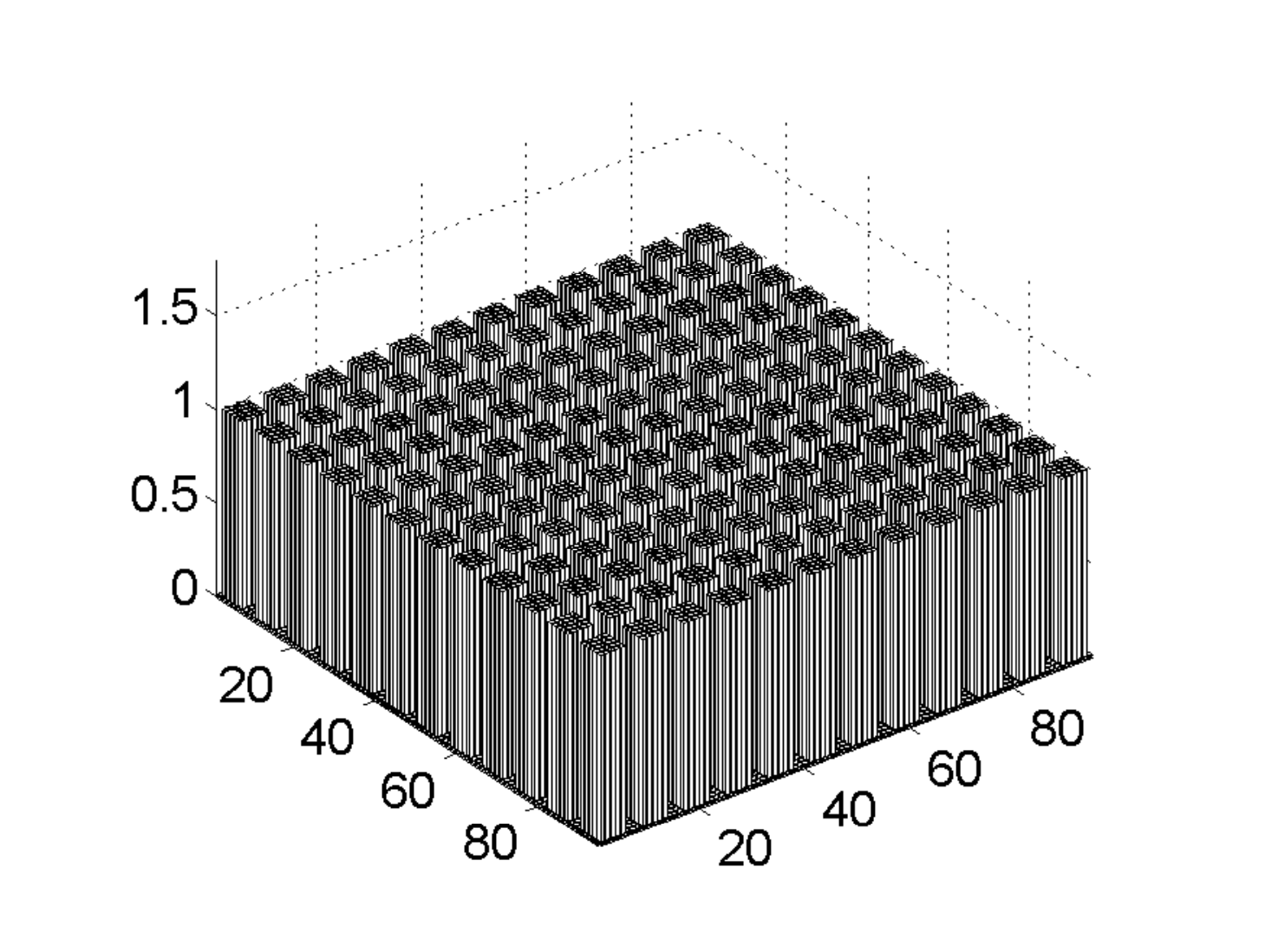}
\includegraphics[width=5.4cm]{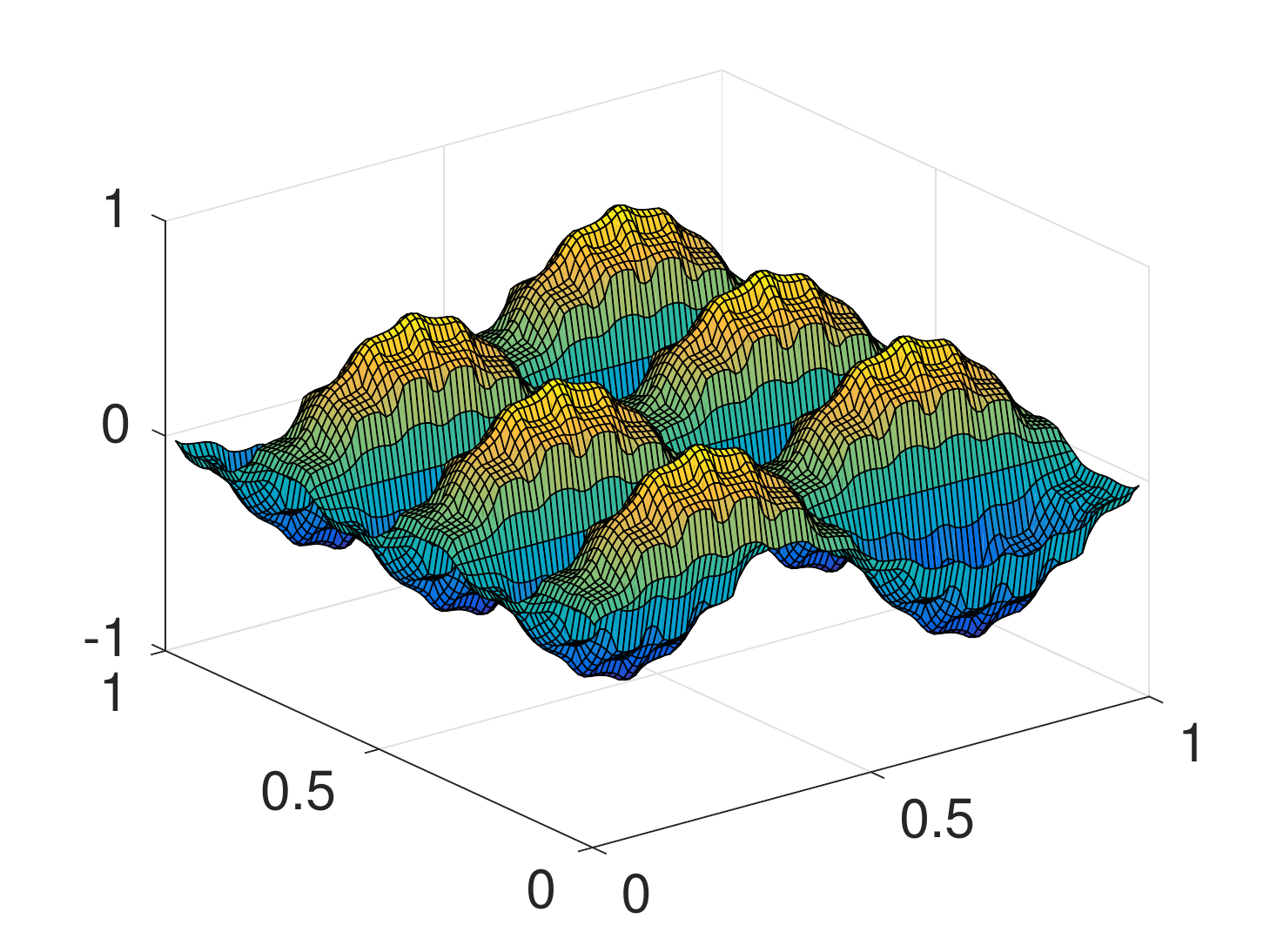}
\includegraphics[width=5.4cm]{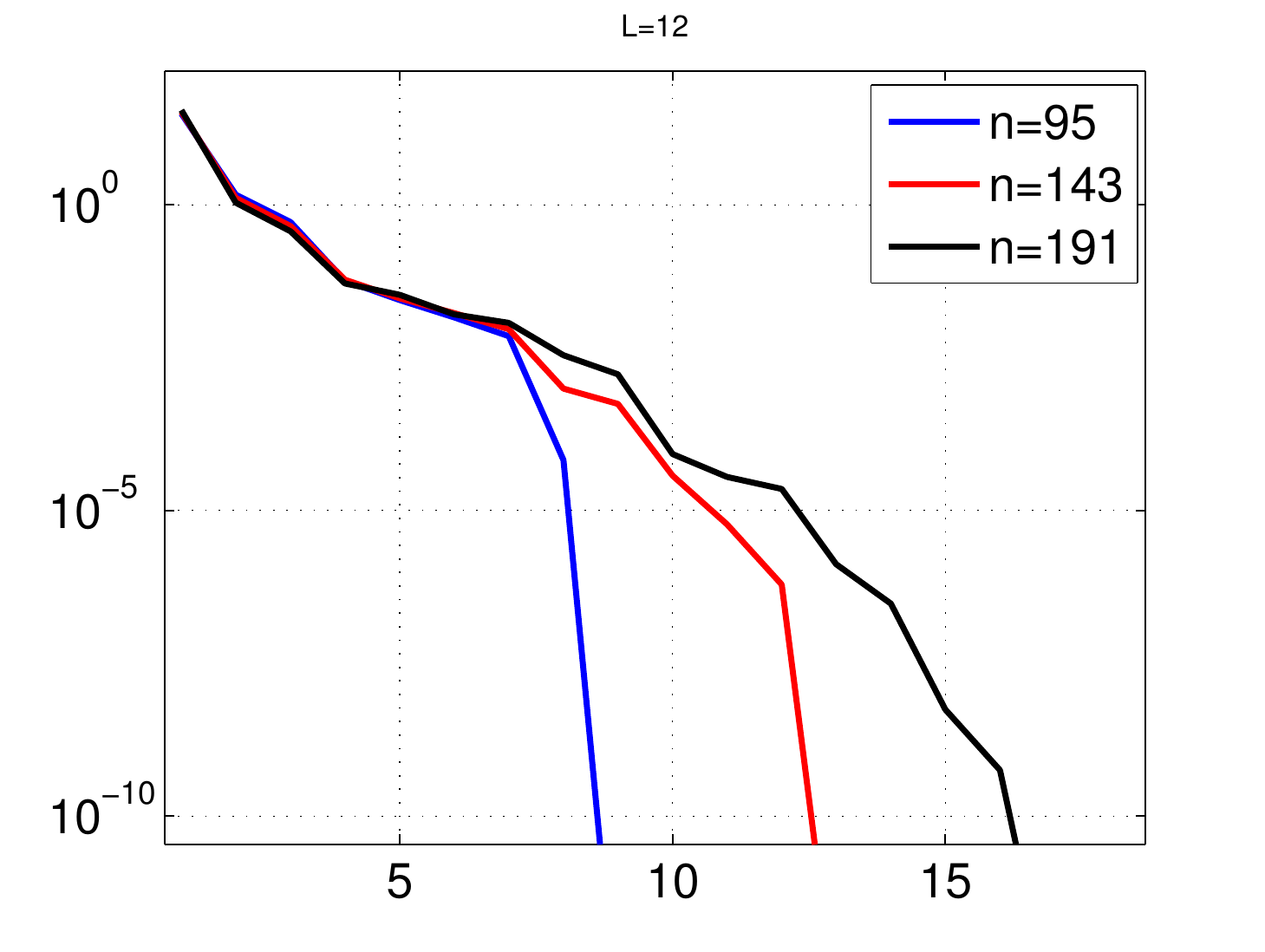}
% %\includegraphics[width=4.0cm]{Coef4StepsSinX-eps-converted-to.pdf}
% %\includegraphics[width=5.0cm]{Sol_L12_n0_8-eps-converted-to.pdf}
% % \includegraphics[width=4.0cm]{CoefSinX3-eps-converted-to.pdf}
% % \includegraphics[width=5.0cm]{Sol_Ranks_L8-eps-converted-to.pdf}
 \caption{\small Rank decomposition of the solution for $12\times 12$ periodic coefficient.}
\label{fig:Rank_12x12}  
\end{figure}

\begin{figure}[htbp]
\centering
% %\includegraphics[width=4.0cm]{UexaMUhomo_Sin_G2e14-eps-converted-to.pdf}
\includegraphics[width=5.4cm]{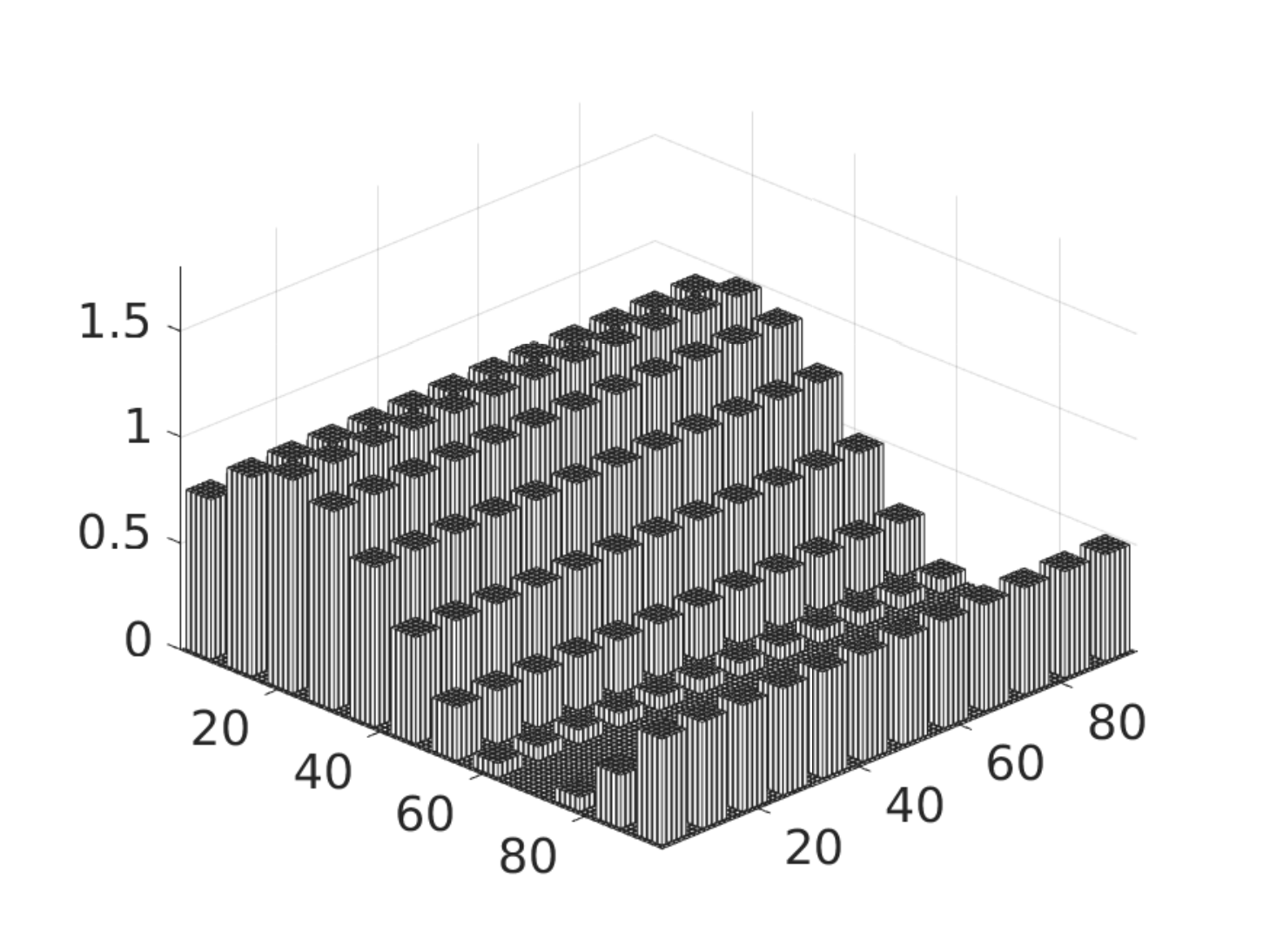}
\includegraphics[width=5.4cm]{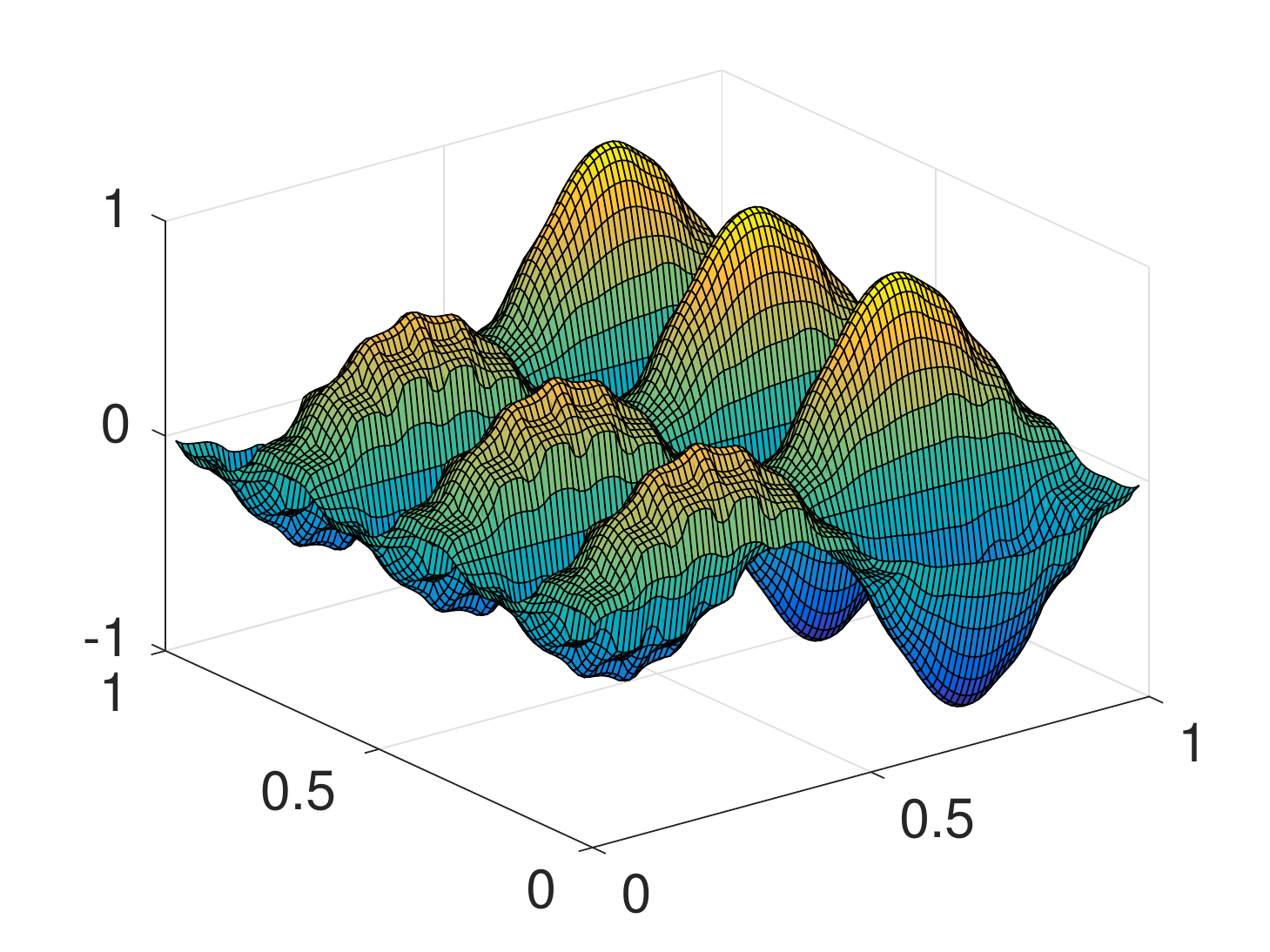}
\includegraphics[width=5.4cm]{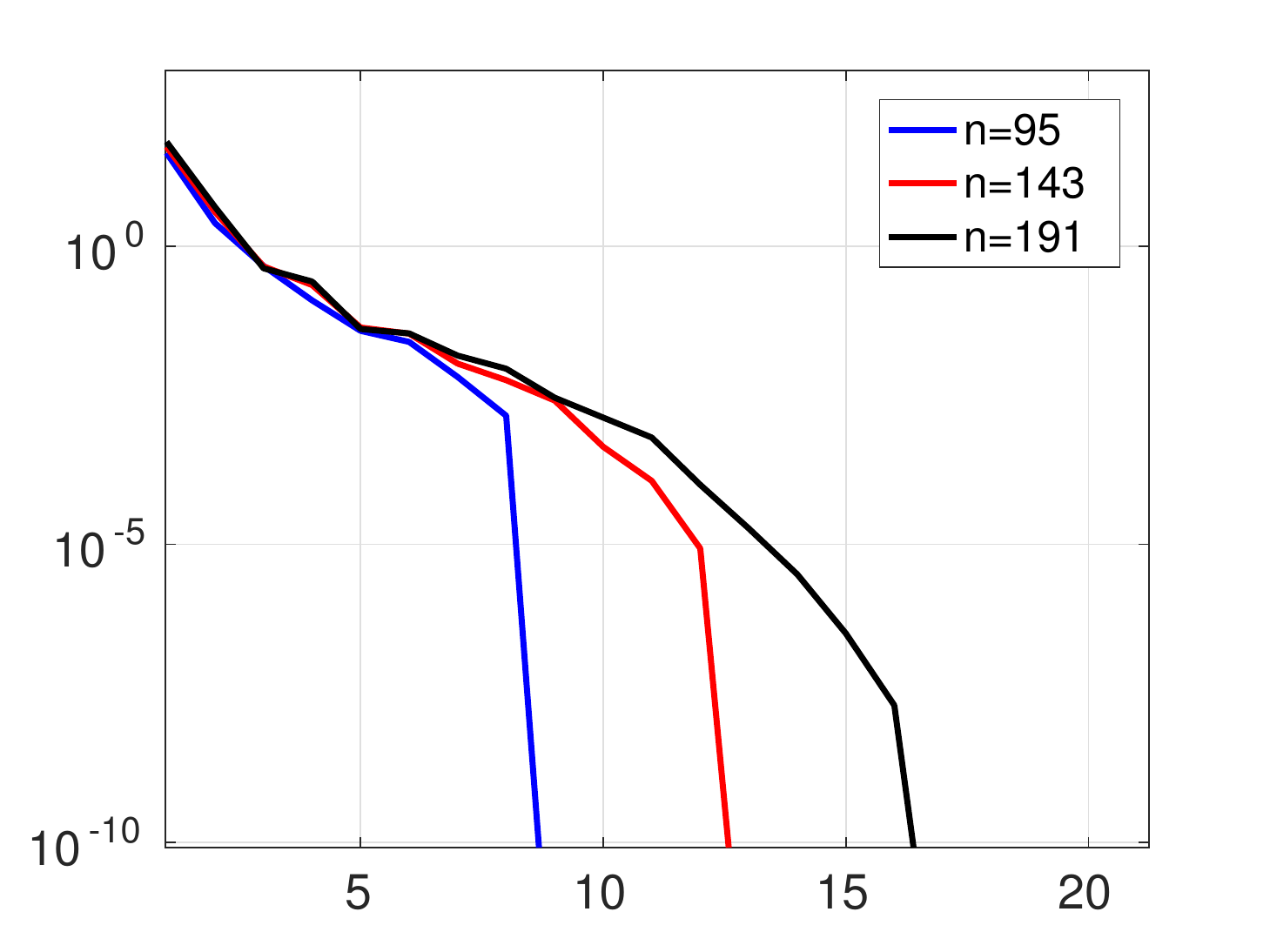}
% %\includegraphics[width=4.0cm]{Coef4StepsSinX-eps-converted-to.pdf}
% %\includegraphics[width=5.0cm]{Sol_L12_n0_8-eps-converted-to.pdf}
% % \includegraphics[width=4.0cm]{CoefSinX3-eps-converted-to.pdf}
% % \includegraphics[width=5.0cm]{Sol_Ranks_L8-eps-converted-to.pdf}
 \caption{\small Rank decomposition of the solution for $12\times 12$ 
 modulated periodic coefficient.}
\label{fig:Rank_12x12_mod}  
\end{figure}

%\color{green}
Figures \ref{fig:Rank_12x12} and \ref{fig:Rank_12x12_mod} demonstrate the singular values 
of the discrete solution on the $n\times n$ grid for $n=95,143, 191$ indicating very moderate
dependence of the $\epsilon$-rank on the grid size $n$. 
As in the case of Figure \ref{fig:Hom2d_Exm1}, in above figures we represent the only oscillating part
of the coefficients and omit the small constant $C>0$.
%\cn

%\color{green}
Further enhancement of the tensor approximation can be based on the application 
of the quantized-TT (QTT) tensor approximation which has been already applied in \cite{BokhSRep:15}
to the 1D equations with quasi-periodic coefficients. 
The power of QTT approximation method is due to the
perfect low rank decompositions applied to the wide class of function-related 
tensors \cite{KhQuant:09}, see \cite{BokhSRep:15} for the more detailed discussion
and a number of numerical examples. 
%\cn

% \cred May be better:\\
% QTT approximation method is highly efficient because a wide class of tensor
% functions interesting in physical applications admit low rank representations  \cite{KhQuant:09}.
% \cn
% Such tensors are obtained by sampling 
% a continuous functions over \cred a \cn uniform (or properly refined) grid.
% \cred
% Frankly I do not understand this sentence. Arbitrary continuous function may
% be so complicated that no low rank representation is possible regardless
% of a grid (take a function with singularity!) 
% \cn

%In our application the QTT tensor approximation applies to the different classes of 
%quasi-periodic equation coefficient. 
%\color{blue}
One can apply QTT approximations to problems with
quasi periodic coefficients, which can be described by 
oscillation with smooth modulation around a constant value, oscillation around 
a given smooth function, or oscillation around piecewise constant function, see Figure 
\ref{fig:1DPeriodStruct2} and examples in \cite{BokhSRep:15}.
%\cn

% \cred
% The structure on the right figure is not clear (too small details). Again minimal values
% seem to be zero!
% \cn
Let the vector ${\bf x}\in \mathbb{C}^N$, $N=2^L$, be obtained by sampling a continuous 
function $f\in C[0,1]$ (or even piecewise smooth functions),
%in C there exists such functions that impossible to present in such simple way!\cn
 on the uniform grid of size $N$. 
 For the following examples of univariate functions the explicit QTT-rank estimates 
 of the corresponding QTT tensor representations
 %\cred which estimates? number  should be given! \cn 
 are valid uniformly  in the vector size $N$, see \cite{KhQuant:09}:\\
%  \cred 
%  you mean uniformly with respect to N? the reader asks why this is important, 
%  where they use this uniformity?
%  \cn
%\begin{itemize}
 (A) $r=1$ for complex exponentials, $f(x)=e^{i \omega x}$, $\omega \in \mathbb{R}$.\\
 (B) $r=2$ for trigonometric functions, $f(x)=\sin {\omega x}, f(x)=\cos {\omega x}$, 
 $\omega \in \mathbb{R}$. \\
 (C) $r\leq m+1$ for polynomials of degree $m$.\\
 (D) For a function $f$ with the QTT-rank $r_0$ modulated by another function $g$ with 
the QTT-rank $r$ (say, step-type function, plain wave, polynomial) 
the QTT rank of a product $f\, g$ is bounded by a multiple of $r$ and $r_0$,
\[
 rank_{QTT}(f g)\leq rank_{QTT}(f )rank_{QTT}(g).
\]
 (E) Furthermore, the following result holds (\cite{VeBoKh:Ewald:14}): 
 QTT rank for the periodic amplification of a 
reference function on a unit cell to a rectangular lattice is of the same order as that 
for the reference function.

%  \begin{proposition} \label{prop:KhQuant_Appr:09}  (\cite{KhQuant:09})
%  Let vector ${\bf x}\in \mathbb{C}^N$, $N=2^L$, be obtained by sampling a continuous function $f\in C[0,1]$ on the uniform grid
%  of size $N$. Then the following QTT-rank estimates are valid uniformly in the vector size $N$:\\
% %\begin{itemize}
% %  (A) $r=1$ for complex exponentials, $f(x)=e^{i \omega x}$, $\omega \in \mathbb{R}$.\\
% %  (B) $r=2$ for trigonometric functions, $f(x)=\sin {\omega x}, f(x)=\cos {\omega x}$, 
% %  $\omega \in \mathbb{R}$. \\
% %  (C) $r\leq m+1$ for polynomials of degree $m$. \\
% % %$r=2$ for Chebyshev polynomials sampled on Chebyshev-Gauss-Lobatto grid.
% (A)
% For a function $f$ with the QTT-rank $r_0$ modulated by another function $g$ with 
% the QTT-rank $r$ (say, step-type function, plain wave, polynomial) 
% the QTT rank of a product $f\, g$ is bounded by a multiple of $r$ and $r_0$,
% \[
%  rank_{QTT}(f g)\leq rank_{QTT}(f )rank_{QTT}(g).
% \]
%  (B) (\cite{VeBoKh:Ewald:14}) QTT rank for the periodic amplification of a 
% reference function on a unit cell to a rectangular lattice is of the same order as that 
% for the reference function.
% %\end{itemize}
% %All estimates are valid uniformly in the vector size $N$. 
% \end{proposition}

The rank of the QTT tensor representation to the 1D Galerkin FEM matrix 
in the case of oscillating coefficients 
was discussed in \cite{DoKazKh_1DSPDE:12,BokhSRep:15}.

\subsection{Numerical test on the rank decomposition of ${\bf u}$}\label{ssec:LowRankSolut}

\begin{figure}[htbp]
\centering
\includegraphics[width=7.0cm]{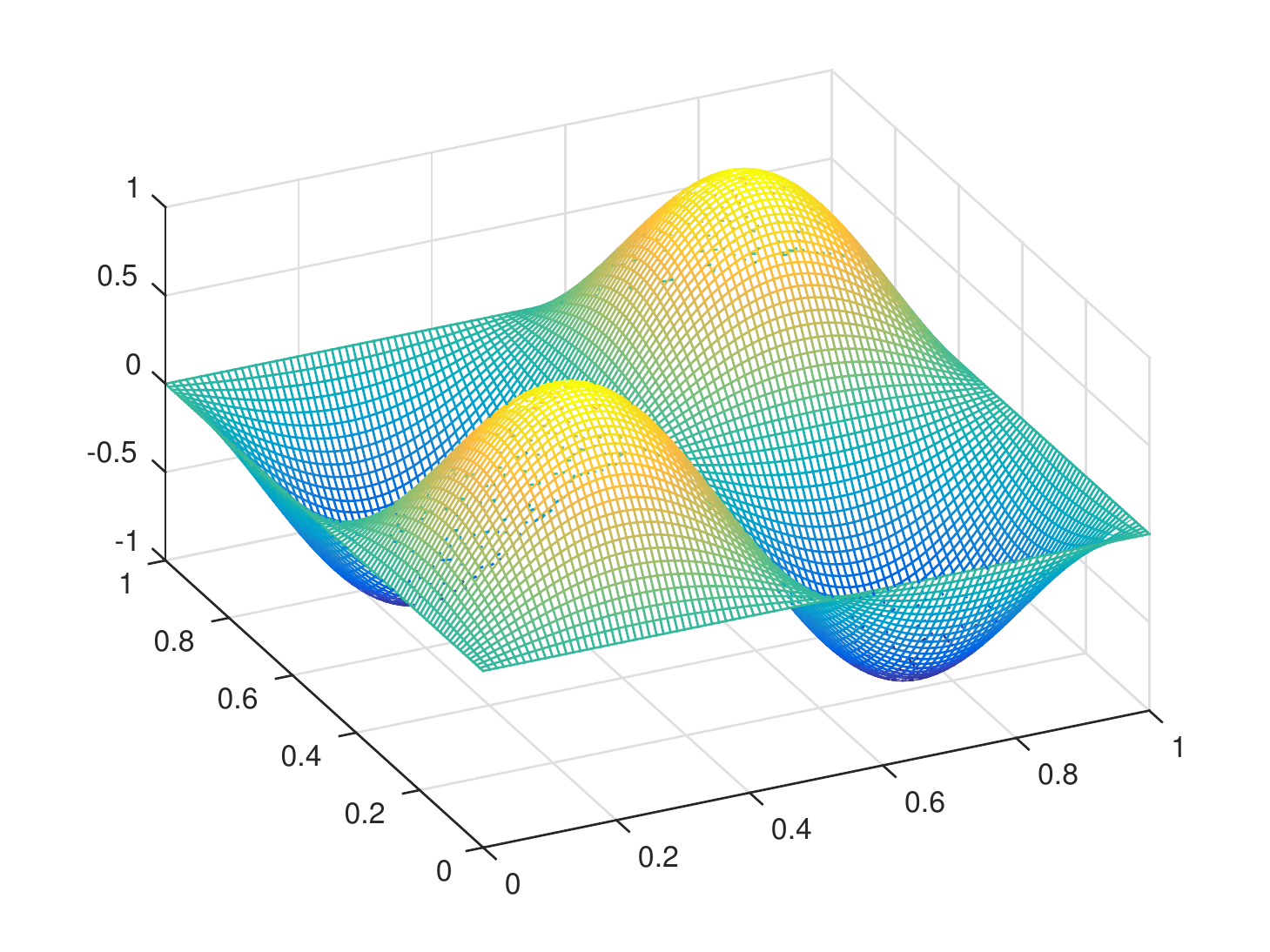}\quad
\includegraphics[width=7.0cm]{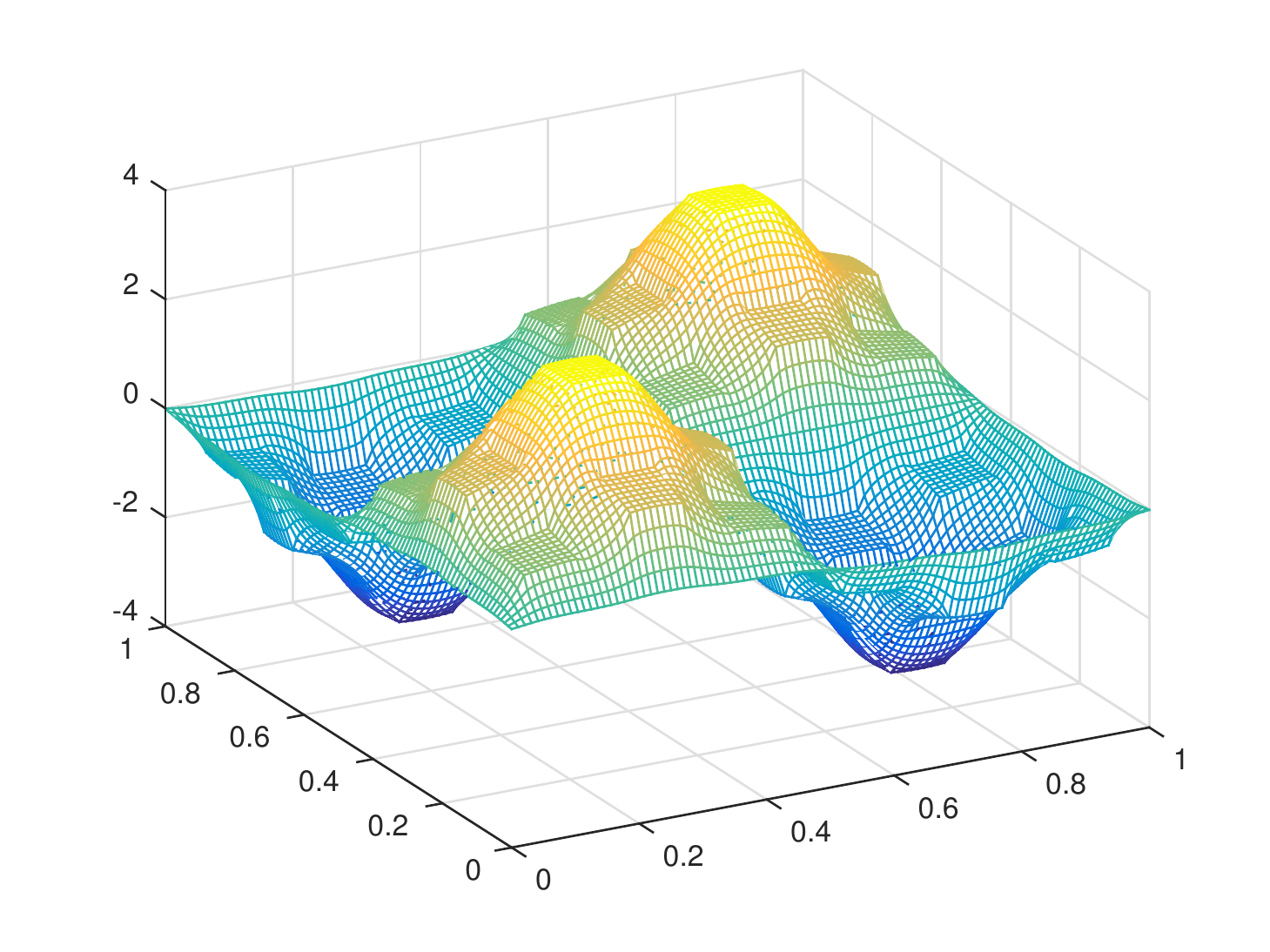}
\caption{\small The right-hand side and solution for periodic oscillating coefficients shown 
in Figure \ref{fig:Hom2d_Exm1}.}
\label{fig:Hom2d_Sol1}  
\end{figure}

\begin{figure}[htbp]
\centering
\includegraphics[width=7.0cm]{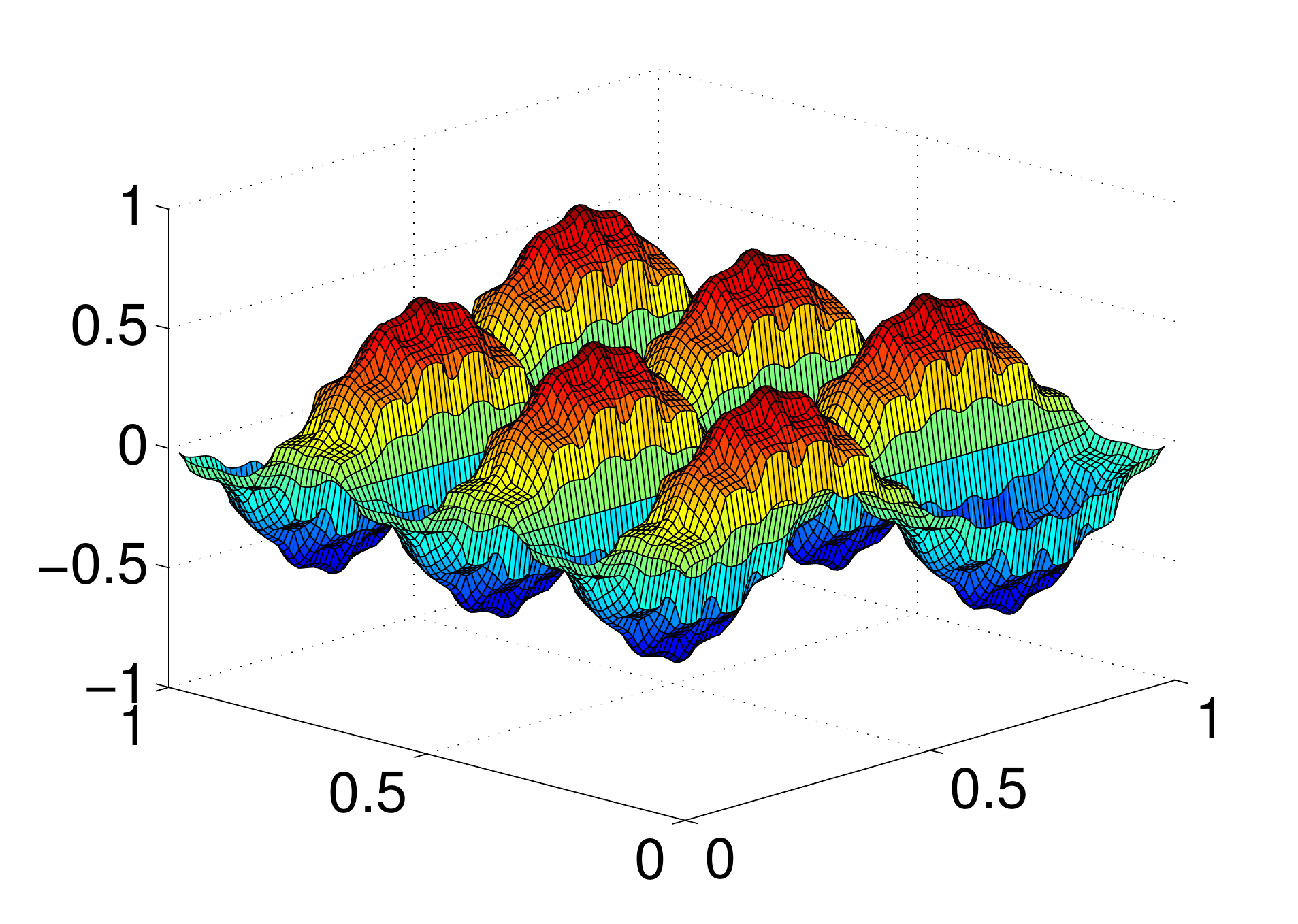}
\includegraphics[width=7.0cm]{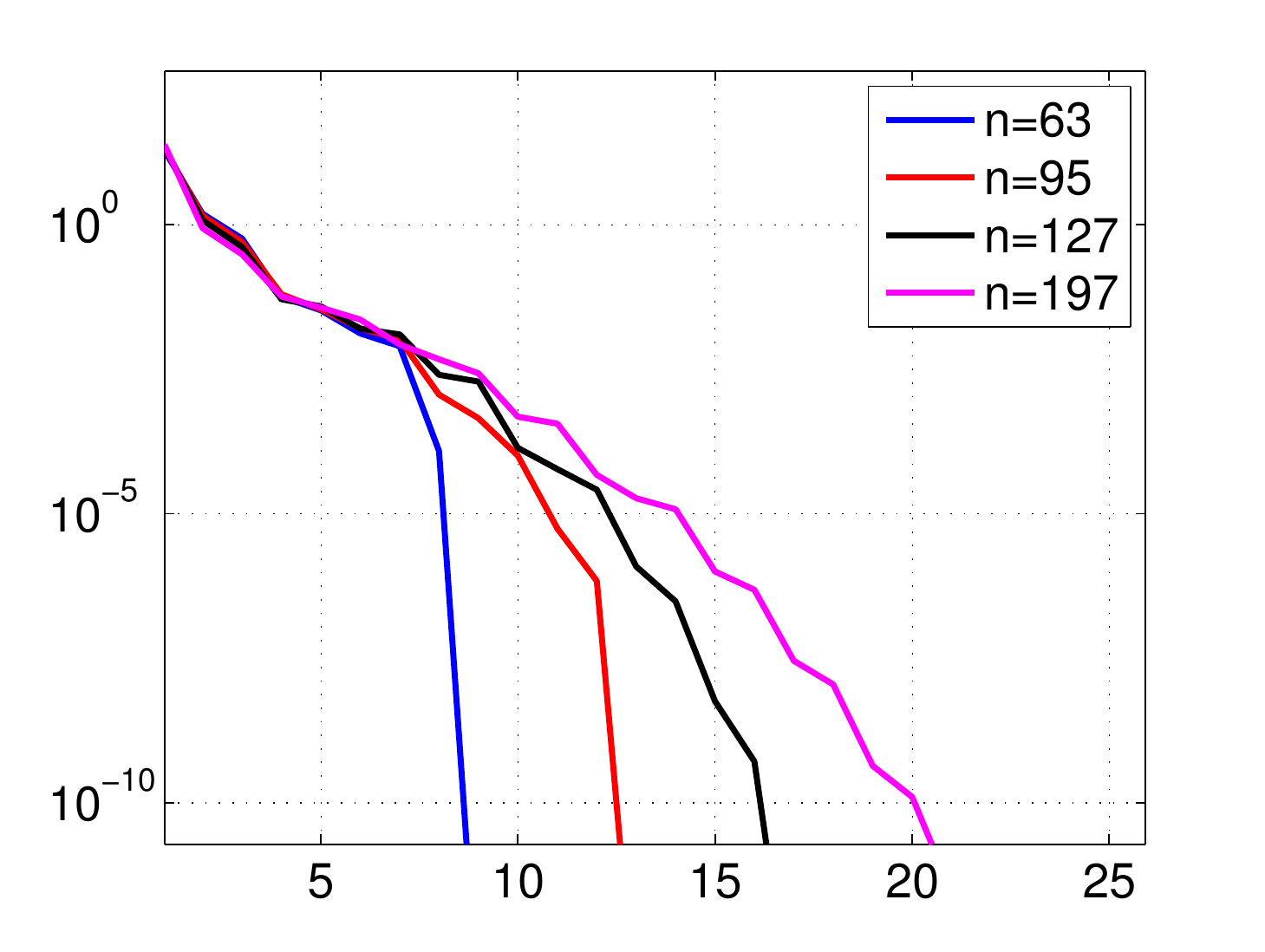}
\caption{\small Accuracy of the rank decomposition of the solution vs. rank parameter 
for $8\times 8$ periodic coefficient and grid size $n\times n$.}
\label{fig:Rank_8x8} 
\end{figure}

Figure \ref{fig:Hom2d_Sol1} represents the right-hand side $f_1(x_1,x_2)$ 
and the respective solution for 
the discretization to equation (\ref{eqn:5.1}) (with the coefficient depicted in Figure 
\ref{fig:Hom2d_Exm1}) on $400 \times 400$-grid, where
\[
 f_1(x_1,x_2) = \sin(2x_1) \sin(2 x_2). %; \quad f_2(x_1,x_2) = \sin(2x_1). 
\]
The PCG solver for the system of equations (\ref{eqn:FEM_syst}) with the 
discrete Laplacian inverse 
as the preconditioner demonstrates robust converges with the rate $q\ll 1$.
Next example demonstrates the rank behavior in the singular value decomposition (SVD) 
of a matrix 
representing the solution vector ${\bf u}\in \mathbb{R}^{n_1\times n_2}$ 
to the equation (\ref{eqn:FEM_syst})
with $12\times 12$ periodic coefficient shown in Figure \ref{fig:Rank_12x12}, left.
Figure \ref{fig:Rank_8x8} represents the rank behavior in the SVD decomposition of 
the solution in the case of $8 \times 8$ periodic coefficient. 

It is worth to observe that comparison of Figures \ref{fig:Rank_12x12} and \ref{fig:Rank_8x8} 
indicates that the exponential decay of the approximation error in the rank 
parameter is stable with respect to the size of $L \times L$ lattice structure of 
the coefficient, i.e. the behavior of the singular values remains almost the 
same for different parameters $\epsilon=1/L$.

Our iterative scheme includes only the matrix-vector multiplication with the 
stiffness matrix $A$ that has the small Kronecker rank $2R$, and the 
action of the preconditioner
defined by the approximate inverse to the Laplacian type matrix. The latter has low 
Kronecker rank of order $R_B=O(|\log \varepsilon |^2)$ as shown above.

Given rank-$1$ vector ${\bf u}={\bf u}_1 \otimes {\bf u}_2$,
the standard property of the Kronecker product matrices 
\[
 {A}{\bf u} = {A}_1{\bf u}_1 \otimes M_2 {\bf u}_2 + M_1 {\bf u}_1 \otimes {A}_2 {\bf u}_2,
\]
indicates that the matrix-vector multiplication enlarges the initial rank by the factor of $2$
and similar with action of preconditioner.
Hence each iterative step should be supplemented with certain rank truncation procedure
which can be implemented adaptively to the chosen approximation threshold or fixed bound
on the rank parameter.

% \cred Comment. I suggest to first finish all the story about our example for $d=2$, put
% examples, comments, pictures, etc. Then as a remark we can extend to higher
% dimensions. DONE BELOW. \cn

\begin{remark}\label{rem:Matr_general_d}
Notice that for $d=3$ the transformed matrix $A$ takes a form
\[ %\ben \label{eqn:Lapl_Kron3}
A = A_1 \otimes I_2\otimes I_3 + I_1 \otimes A_2 \otimes I_3 + I_1 \otimes I_2\otimes A_3,
\]  % \een
and it obeys the $d$-term Kronecker sum representation in the general.
Hence, in the general case of $d\geq 2$ and $R\geq 1$ the Kronecker rank of the 
matrix ${A}$ is given by
\[
 \mbox{rank}_{Kron} ({A})= d \, R.
\]
\end{remark}

\section{Conclusions}\label{sec:Conclusions}

We present a preconditioned iteration method for solving an elliptic type boundary value problem
in $\mathbb{R}^d$
with the operator generated by a quasi--periodic structure
with rapidly changing coefficients characterized by a small length parameter $\epsilon$. 
We use tensor product FEM discretization
that allows to approximate the stiffness matrix $A$ in the form of low-rank Kronecker sum.
The preconditioner $\mathcal{A}_0$
is constructed based on certain averaging (homogenization) procedure of the initial
equation coefficients such that inversion of $\mathcal{A}_0$ is much simpler 
than inversion of $\mathcal{A}$. 
We prove contraction of the iteration method and establish explicit estimates of the contraction
factor $q<1$.
For typical quasi--periodic structures we deduce fully computable two--sided a 
posteriori estimates which are able to control numerical solutions on any iteration.

We apply the tensor-structured approximation which is especially efficient
if the equation coefficients admit low rank representations and algebraic operations
are performed in tensor structured formats.
Under moderate assumptions the storage and solution complexity
of our approach depends only weakly (merely linear-logarithmically) 
on the frequency parameter $1/\epsilon$. 
Numerical tests demonstrate that the FEM solution allows the accurate 
low rank separable approximation which is the basic prerequisite
for application of the tensor numerical methods to the problems of geometric homogenization.

The approach allows further enhancement based on the quantized-TT (QTT)
tensor approximation which is the topic for future research work.
Another direction is related to fully tensor structured implementation of the 
computable two--sided a posteriori error estimates.

\vspace{0.3cm}
{\bf Acknowledgements.} 
SR appreciates the support provided by the Max-Planck Institute
for Mathematics in the Sciences (Leipzig, Germany) during his scientific visit in 2016.
The authors are thankful to Dr. V. Khoromskaia (MPI MIS, Leipzig) 
for  the numerical experiments.

% \section{Numerical examples}
% \vskip3cm
% {TO BE ADDED}
% \vskip3cm
% {............................................}
% 

%%%%%%%%%%%%%%%%%%%%%%%%%%%%%%%%%%%

%%%%%%%%%%%%%%%%%%%%%%%%%%%%%

% \bibitem{Vidal2003} G. Vidal.
%   \emph{Efficient classical simulation of slightly entangled quantum
% computations.} Phys. Rev. Lett. 91(14), 2003, 147902-1 147902-4.

% \bibitem {OsTy_TT:09} I.V. Oseledets, and E.E. Tyrtyshnikov.
%   \emph{Breaking the Curse of Dimensionality, or How to Use SVD in
%     Many Dimensions.} SIAM J. Sci. Comp., 31 (2009), 3744-3759.
% \bibitem{Cirac_TC:04} F. Verstraete, D. Porras, and J.I. Cirac.
% \emph{DMRG and periodic boundary conditions: A quantum information perspective. }
% Phys. Rev. Lett., 93(22): 227205, Nov. 2004.

% \bibitem{Vidal2003} G. Vidal.
%   \emph{Efficient classical simulation of slightly entangled quantum
% computations.} Phys. Rev. Lett. 91(14), 2003, 147902-1 147902-4.
%%%%%%%%%%%%%%%%%%%%%%%%%%%%%%
%%%%%%%%%%%%%%%%%%%%%%%%%%%%%%
%%%%%%%%%%%%%%%%%%%%%%%%%%%%%%


\begin{thebibliography}{99}                                                                                              

\bibitem {Bakhvalov} Bakhvalov, N. S., Panasenko, G. Homogenisation:
Averaging Processes In Periodic Media: Mathematical Problems In The Mechanics
Of Composite Materials. Springer, 1989.

\bibitem{BeKh3_Prot1:16} P. Benner, V. Khoromskaia and B. N. Khoromskij.
\emph{Range-separated  tensor formats for numerical modeling of many-particle 
interaction potentials}.
E-preprint, http://arxiv.org/abs/1606.09218, 2016.

\bibitem {Bensoussan} Bensoussan, A., Lions, J.-L., Papanicolaou, G. (1978):
Asymptotic analysis for periodic structures. Amsterdam: North-Holland

\bibitem{BrennSc} S. Brenner  and R. Scott. \emph{The mathematical theory of finite
element methods.} Springer, 1994.

\bibitem{Davis} J. P. Davis. Circulant matrices. New York. John Wiley \& Sons, 1979.


%
%
%\bibitem {Chipot} Chipot, M. (2009): Elliptic Equations: An Introductory
%Course. Birkh\"{a}user Verlag AG

%\bibitem {Cioranescu}
%Cioranescu, D., Donato, P. (1999): An Introduction to
%Homogenization. Oxford Lecture Series in Mathematics and its Applications. Bd.
%17. Oxford University Press

%\bibitem{25} D{\"o}rfler, W., Rumpf, M. (1998): 
%An adaptive strategy for elliptic problems including a posteriori controlled boundary approximation. 
%Math.Comput., \textbf{67}, 224, 1361-1382


\bibitem {Jikov} Jikov, V.V., Kozlov, S.M., Oleinik, O.A. (1994):
Homogenization of differential operators and integral functionals. Berlin: Springer
\bibitem {Friedman}Friedman, A. (1976): Partial Differential Equations. R. E.
Krieger Pub. Co., Huntington, NY

\bibitem{GaHaKh4:02}
I.P. Gavrilyuk, W. Hackbusch and B.N. Khoromskij. {\it  
  Hierarchical Tensor-Product Approximation to the Inverse and
  Related Operators in High-Dimensional Elliptic Problems. }
Computing {\bf 74} (2005), 131-157. 

\bibitem{GlOtto:15}
Antoine Gloria and Felix Otto. 
Quantitative estimates on the periodic approximation of the corrector in stochastic homogenization
In: ESAIM / Proceedings, 48 (2015), p. 80-97.\\
MIS-Preprint  12/2015, DOI: 10.1051/proc/201448003.

\bibitem{GlLiTr}
R. Glowinski, J.-L. Lions, R. Tr\'emolier\'es.
Analyse num\'erique des in\'equations variationnelles.
Dunod, Paris, 1976. 

\bibitem{KantorovichKrylov}
Kantorovich L. V. and Krylov V. L., Approximate
Methods of Higher Analysis. Interscience, New York, 1958.

\bibitem{KaOsRaSch:16}
V. Kazeev, I. Oseledets, M. Rakhuba, and Ch. Schwab.
\emph{QTT-finite-element approximation for multiscale
problems I: model problems in one dimension.}
Adv. Comput. Math., 2016. DOI: 10.1007/s10444-016-9491-y.


\bibitem{KaRaSch:11}
V. Kazeev,  O. Reichmann, and Ch. Schwab.
\emph{Low-rank tensor structure of linear diffusion operators in the TT and QTT formats.}
Linear Algebra and its Applications, v. 438(11), 2013, 4204-4221.

\bibitem {DoKazKh_1DSPDE:12} S. Dolgov, V. Kazeev, and B.N. Khoromskij.
\emph{The tensor-structured solution of one-dimensional elliptic differential equations 
with high-dimensional parameters.}
Preprint 51/2012, MPI MiS, Leipzig 2012.

\bibitem{VeBoKh:Ewald:14} V. Khoromskaia and B. N. Khoromskij.
\emph{Grid-based lattice summation of electrostatic potentials by assembled rank-structured 
tensor approximation.} Comp. Phys. Commun., {\bf 185} (12), 2014, pp. 3162-3174. 
%DOI: 10.1016/j.cpc.2014.08.015. 
%Preprint arXiv:1405.2270. 

\bibitem{VeKhorCorePeriod:14} V. Khoromskaia, and B.N. Khoromskij.
\emph{Tensor Approach to Linearized Hartree-Fock Equation  for Lattice-type and Periodic Systems.}
%Preprint 62/2014, MPI MiS, Leipzig 2014.
E-preprint arXiv:1408.3839, 2014.


\bibitem{VeKhorEwTuck_NLLA:15} V. Khoromskaia and B.N. Khoromskij.
\emph{Fast tensor method for summation of long-range potentials on 3D lattices with defects}.
Numerical Linear Algebra with Applications, 2016, v. 23: 249-271. 
%DOI: 10.1002/nla.2023.
% %Preprint 65/2015, MPI MiS, Leipzig 2015.

\bibitem{VeKhorTromsoe:15} V. Khoromskaia and B.N. Khoromskij.
\emph{Tensor numerical methods in quantum chemistry: from Hartree-Fock to excitation energies.}
Phys. Chem. Chem. Phys., 17:31491 - 31509, 2015. % DOI:10.1039/c5cp01215e.
%E-preprint arXiv:1504.06289, 2015.
%Preprint 19/2015, MPI MiS, Leipzig 2014 (submitted).

% \bibitem{KhKhOt_StHomo:16} V. Khoromskaia, B. N. Khoromskij and F. Otto.
% \emph{A numerical primer in 2D stochastic homogenization.}
% Manuscript in preparation, 2016.

\bibitem{KhorCA:09} B.N. Khoromskij.
  \emph{Tensor-Structured Preconditioners and Approximate Inverse
of Elliptic Operators in  $\mathbb{R}^d $.} J. Constr. Approx. {\bf 30} (2009) 599-620.

\bibitem{KhQuant:09} B.N. Khoromskij.
\emph{$O(d\log N)$-Quantics Approximation of $N$-$d$ Tensors in High-Dimensional
Numerical Modeling.}
%Preprint 55/2009 MPI MiS, Leipzig 2009.\\
Constr. Approx. 34  (2011) 257--280.

\bibitem{KhorSurv:10}   B.N. Khoromskij. \emph{Tensors-structured
Numerical Methods in Scientific Computing: Survey on Recent Advances.}
Chemometr. Intell. Lab. Syst. 110 (2012), 1-19.
%DOI: 10.1016/j.chemolab.2011.09.001.}

\bibitem{KhWi:B}
B.N. Khoromskij and G. Wittum. {\it Numerical Solution of Elliptic
Differential Equations by Reduction to the Interface}.
Research monograph,  LNCSE, No. 36, Springer-Verlag, 2004.

% \bibitem{KhorSurv:14} Boris N. Khoromskij.
% \emph{Tensor Numerical Methods for High-dimensional PDEs: Basic Theory and
%   Initial Applications}. ESAIM: Proceedings and Surveys,  
% Eds. N. Champagnat, T. Leli{\'e}vre, A. Nouy. January 2015, Vol. 48, p. 1-28.

\bibitem{BokhSRep:15}  B.N. Khoromskij and S. Repin.
\emph{A fast iteration method for solving elliptic problems with quasiperiodic coefficients}. 
Russ. J. Numer. Anal. Math. Modelling 2015; 30 (6):329-344.  
E-preprint arXiv:1510.00284, 2015.

\bibitem{KhSautVeit:11} B.N. Khoromskij, S. Sauter, and A. Veit. 
\emph{Fast Quadrature 
Techniques for Retarded Potentials Based on TT/QTT Tensor Approximation.} 
Comp. Meth. in Applied Math., v.11 (2011), No. 3, 342 - 362.
%Preprint 44/2011, MPI MiS, Leipzig 2011.


\bibitem{LiSt1967}
J.-L. Lions and G. Stampacchia.
 Variational inequalities. Comm. Pure Appl. Math. 20 1967 493--519. 

 \bibitem{OselDolg:12}
 Ivan V Oseledets, and  S.V. Dolgov.
 \emph{Solution of linear systems and matrix inversion in the TT-format. }
 SIAM Journal on Scientific Computing, v. 34(5), 2012, A2718-A2739.
 
\bibitem{MaNeRe}
O. Mali, P. Neittaanmaki, S. Repin.
Accuracy verification methods. Theory and algorithms.
Springer, 2014

\bibitem{MarShai} G. I. Marchuk and V. V. Shaidurov.
\emph{Difference methods and their extrapolations}.
Applications of Mathematics, New York: Springer, 1983.


\bibitem{NeRe}
P. Neittaanmaki and S. Repin.
Reliable methods for computer simulation.
Error control and a posteriori estimates. Elsevier, 2004.

\bibitem{Ostrowski}
A. Ostrowski. Les estimations des erreurs a posteriori dans les proc\'ed\'es it\'eratifs,
C. R. Acad. Sci, Paris, S\'er. A–B 275 (1972), pp. A275–A278.


\bibitem{Re2000} S. Repin. \emph{A posteriori error estimation for variational
problems with uniformly convex functionals,} { \em Math. Comput.},
69(2000),  230,   481--500.


\bibitem{ReGruyter}
S. Repin. \emph{A Posteriori Estimates for Partial Differential Equations}.
Walter de Gruyter, Berlin, 2008.

\bibitem{ReSaSa1}
S. Repin, T. Samrowski, and S. Sauter.
A posteriori error majorants of the modeling errors for elliptic homogenization problems.
\emph{C. R. Math. Acad. Sci. Paris} 351 (2013), no. 23-24, 877-882

\bibitem{ReSaSa2}
S. Repin, T. Samrowski, and S. Sauter.
Combined a posteriori modeling-discretization error estimate for elliptic problems with complicated interfaces.
\emph{ESAIM Math. Model. Numer. Anal.}, 46 (2012), no. 6, 1389-1405.

\bibitem{ReSaSm}
S. Repin, S. Sauter, and A. Smolianski. A posteriori estimation of dimension reduction errors for elliptic problems on thin domains.
SIAM J. Numer. Anal. 42 (2004), no. 4, 1435--1451.

\bibitem{Scholl:11} U. Schollw\"ock. \emph{The density-matrix renormalization group in the age of
matrix product states}, Ann.Phys. 326 (1) (2011) 96-192.

% \bibitem{Cirac_TC:04} F. Verstraete, D. Porras, and J.I. Cirac.
% \emph{DMRG and periodic boundary conditions: A quantum information perspective. }
% Phys. Rev. Lett., 93(22): 227205, Nov. 2004.
% 
% \bibitem {White:93} S.R. White.
% \emph{Density-matrix algorithms for quantum renormalization groups.}
% Phys. Rev. B, v. 48(14), 1993, 10345-10356.

\bibitem{Zeidler}
E. Zeidler. Nonlinear functional analysis and its applications.
I. Fixed-point theorems, Springer-Verlag, New York, 1986.

  \end{thebibliography}
\end{document}